\input amstex
\input amsppt.sty
\loadbold
\magnification=1200
\def\la{{\lambda}}
\def\bla{{\boldsymbol \lambda}}

\def\sB{{\Cal B}}

\def\sD{{\Cal D}}

\def\sH{{\Cal H}}

\def\sO{{\Cal O}}

\def\bH{{\bold H}}

\def\bS{{\bold S}}

\def\al{\alpha} 
 
\def\de{\delta}

\def\la{\lambda}

\def\La{\Lambda}

\def\le{\leqslant}
\def\ge{\geqslant}
\def\fS{\frak S}
\def\End{\text{\rm End}}

\def\dim{\text{\rm dim\,}}

\def\hmod{\text{-\bf mod}}

\def\lecq{\preccurlyeq}

\def\recq{\succcurlyeq}
\def\rec{\succ}

\def\wt{{\text{\rm wt}}}

\def\bka{{\boldkey a}}
\def\bkb{{\boldkey b}}
\def\bkc{{\boldkey c}}
\def\op{{\text{\rm op}}}
\def\cns{{\text{\rm cns}}}
\def\proclaim#1{\vskip.3cm\noindent{\bf#1.}\quad\it}
\def\endproclaim{\par\vskip.3cm\rm}
\def\tpi{{\tilde\pi}}

\topmatter
\title Ariki-Koike  algebras with semisimple bottoms  \endtitle
\rightheadtext{Ariki-Koike  algebras}
\thanks The research was carried out while the second author
was visiting the University of New South Wales.
Both authors gratefully acknowledge the support received from 
the Australian Research
Council. The second author is partially
supported by the National Natural Science Foundation 
in China. He wishes to thank the University of New South Wales for its
hospitality during the writing of the paper.\endthanks
\author Jie Du  and Hebing Rui  \endauthor
\affil {\eightpoint
School of Mathematics, University of New South Wales\\  Sydney,
 2052, Australia \\ 
 Department of Mathematics,   University of Shanghai for Science \&
Technology\\
Shanghai, 200093, P.R.China}\endaffil 
%\date 30 July, 1998\enddate
\email  jied\@maths.unsw.edu.au and hbruik\@online.sh.cn\endemail
%\vsize22.5truecm
%\baselineskip15pt
\endtopmatter
\vsize24truecm

Consider the cyclic group $C_m$ of order $m$ and the wreath product 
$W=W_m^r=C_m\wr \fS_r$, where $\fS_r$ is the symmetric group  on $r$ letters.
Then the direct product  $C=C_m\times\cdots\times C_m$ of $r$ copies
$C_m$ is normal in $W$ and sits at the ``bottom'' of $W$.
Let $F$ be a splitting field of $C$ in which
$m$ is not zero. Then the ``bottom''
$FC$ of $FW$ is semisimple and simple $C$-modules
$N_{\boldkey i}$ are indexed by the set $I(m,r)$ of all
$r$-tuples $\boldkey i=(i_1,\cdots,i_r)$ with $1\le i_j\le m$.   
Moreover, there is a central primitive idempotent decomposition 
$1=\sum_{\boldkey i\in I(m,r)}e_{\boldkey i}$ such that 
$N_{\boldkey i}\cong e_{\boldkey i}FC$.
The symmetric group $\fS_r$ acts on $I(m,r)$ by place permutation,
and the set of $W$-orbits is identified with the set $\La(m,r)$ of 
compositions of $r$ with $m$ parts. Thus $\boldkey i$ is in the 
orbit $\la\in\La(m,r)$, denoted $\wt(\boldkey i)=\la$, if
$\la_j=\#\{i_k:i_k=j\}$ for all $j$. Notice that, if $\wt(\boldkey i)=
\wt(\boldkey j)$, then we have isomorphism of induced modules:
$N_{\boldkey i}\uparrow^W\cong N_{\boldkey j}\uparrow^W$.
Therefore, putting $N_\la=N_{\boldkey i}$ and
$e_\la=e_{\boldkey i}$ if $\wt(\boldkey i)=\la$,
we have right $W$-module isomorphism
  $$FW\cong\bigoplus_{\la\in\La(m,r)}(N_\la\otimes_{FC}FW)^{\oplus d_\la}
\cong\bigoplus_{\la\in\La(m,r)}(e_\la FW)^{\oplus d_\la}.$$
Standard results  will give a Morita 
equivalence between the categories of
$FW$-modules and $e FWe$-modules, where $e=\sum_{\la\in\La(m,r)}e_\la$.
Since $e_\la FW e_\mu\cong \de_{\la\mu}F\fS_\la$, 
where $\fS_\la$ is the Young subgroup
corresponding to $\la$, we have Morita equivalence
$$FW\hmod\overset{\text{Morita}}\to\sim (\oplus_{\la\in\La(m,r)}F\fS_\la)\hmod.
\tag1$$
On the other hand, for the Hecke algebra of type $B$ (i.e.,
the Hecke algebra associated to the
group $W_2^r$),
Dipper and James established a Morita equivalence analogous to (1)  
in \cite{DJ2}.

This paper is going to  generalize these results to the 
Ariki-Koike algebra $\bH=\bH_m^r$, an Iwahori-Hecke type algebra 
associated with $W_m^r$ (see \cite{AK}). 
A major difficulty here is the non-existence
of a subalgebra based on the bottom $C$, comparing with the classical case,
and also, the group $W$ is no longer a Coxeter group.
However, since the semi-simplicity of $FC$ is simply equivalent to
the condition that the order $|C|$ of the bottom $C$
is non-zero in  $F$, our strategy here is to find
a $q$-analogue of the order of $C$, which is
called the Poincar\'e polynomial of $C$, and with the invertibility
of such a polynomial, to look for those idempotents $e_\la$. 
Thus, we eventually establish a $q$-analogue of the Morita equivalence (1)
above. A by-product of our results is the introduction of the 
Poincar\'e polynomial $d_W$ of the complex reflection group $W$.
We shall see that the semi-simplicity of $\bH$ over a field $F$
is equivalent to $d_W\neq0$ in $F$.

We organize the paper as follows.
In \S1, we introduce the poset $\La[m,r]$, which is isomorphic to
$\La(m,r)$ and discuss some combinatorics related to
symmetric groups. In \S2, a useful lemma (2.8)  related to the
poset structure on $\La[m,r]$ is proved. 
Candidates of those idempotents $e_\la$ are constructed in \S3.
The main results are presented in \S4, where we prove that the
invertibility of the `Poincar\'e polynomial' $f_{m,r}$
is  a necessary and sufficient condition for
the existence of those idempotents $e_\la$, and the Morita equivalence
is established. 
Finally, in \S5, we lift the Morita equivalence
to the endomorphism algebra level. Two by-products for $\bH_F$ 
over a field $F$ in which $f_{m,r}$ is nonzero
are the classification of simple modules and the criterion of semisimplicity. 

The main results of the paper have been announced by the first author at the
``Symposium on Modular representations of finite groups'', Charlottesville,
Virginia,
May 1998, and at the ``International Conference on Representation Theory,''
 Shanghai, June-July, 1998.  
At the Virginia conference, R. Dipper announced some Morita theorems
for Ariki-Koike algebras
joint with A. Mathas with quite different treatment. 
For example, our method works over the set $\La(m,r)$ of compositions with
$m$ parts, while, in their method, they first treat the case where compositions
have 2 parts. Thus,
they annouced a Morita equivalence between an Ariki-Koike
algebra and a tensor product of two smaller such algebras.
To obtain our result, they have to break two-part compositions further down. 

Throughout, $R$ denotes a commutative ring with identity 1.

\subhead 1. The poset $\La[m,r]$ \endsubhead
Let $r$ be a non-negative integer.
A {\it composition} $\la$ of $r$ with $m>0$ parts is a sequence
$(\la_1, \cdots, \la_m)$ of nonnegative integers  such that 
$|\la|=\sum_{i=1}^m \la_i=r$, and $\la$ is called a 
{\it partition} if the sequence is  weakly  decreasing.
Let $\La(m, r)$ (resp.
$\La(m, r)^+$) be the set of compositions (resp. partitions) of $r$ with
$m$-parts. 
 With the usual
dominance order $\trianglelefteq$, both $\La(m,r)$ and $\La(m,r)^+$
are posets.

For notational convenience, we shall use another poset $\La[m,r]$.
For any $\la\in\La(m,r)$,
let $[\la]= [a_0, a_1,\cdots, a_m]$ where $a_0=0$ and $a_i=\la_1+\cdots+\la_i$
for all $i$ and put $\La[m, r]=\{[\la]:\la\in\La(m,r)\}$.
The following results are almost obvious.

\proclaim{(1.1) Lemma} (a) Alternatively, we have
$$
\La[m, r]=\{[a_0, a_1,\cdots, a_m] :  0=a_0\le a_1\le\cdots\le  
a_m=r, a_i\in \Bbb Z, \forall i \}.
$$

(b) For any $\boldkey a, \boldkey b\in \La[m, r]$, define $\boldkey a\lecq \boldkey b$ by 
setting $a_i\le b_i$ for every $i$ with $1\le i\le m$.
Then $\La[ m, r]$ is a poset with partial ordering $\lecq$.
Moreover, the map $\Xi:\La(m,r)\to\La[m,r]$ defined by
$\Xi(\la)=[\la]$ is an isomorphism between posets
$(\La(m,r),\trianglelefteq)$ and $(\La[m,r],\lecq)$.
In particular, we have for all $\la,\mu\in\La(m,r)$
$$\la\trianglelefteq\mu\text{ if and only if } [\la]\lecq[\mu].$$

(c) If $\Theta=\Xi^{-1}$ is the inverse map of $\Xi$, then, for  
$\boldkey a=[a_i]\in\La[m,r]$,
$$\Theta(\boldkey a)=(a_1-a_0, a_2-a_1, 
\cdots, a_m-a_{m-1}).$$
\endproclaim
 
If $\bka\lecq\bkb$ and $\bka\neq\bkb$, we write $\bka\prec\bkb$.

\definition{(1.2) Notation} For any $\bka=[a_0,\cdots,a_m]\in\La[m,r]$, let
$i$ (resp. $ j$ ) be the  minimal index such that
$a_i\neq 0$ (resp. $a_j=r$) and define 

(a) $\boldkey a'=[0, r-a_{m-1}, \cdots, r-a_1,r]$;

(b) $\boldkey a_\vdash =[0, \cdots,  0, a_i-1, \cdots,  a_m-1]$;

(c) $\boldkey a_\dashv=[a_0, a_1, \cdots, a_{j-1}, r-1, \cdots, r-1].$

\noindent
The notations we choose here are symmetric: 
if $\la=\Theta(\bka)$, then
 $i$ (resp. $ j$ ) is the  minimal (resp. maximal) index
 with $\la_i\neq 0$ (resp. $\la_j\not=0$), and we have 
$\boldkey a'=[(\la_m,\cdots,\la_1)]=[\la^\circ]$,
$\boldkey a_\vdash=[(0,\cdots,0,\la_i-1,\la_{i+1},\cdots,\la_m)]
=[\la_\vdash]$, and
$\boldkey a_\dashv=[(\la_1,\cdots,\la_{j-1},\la_j-1, 0, \cdots, 0)]
=[\la_\dashv]$.
Since  $(\la_\vdash)^\circ=(\la^\circ)_\dashv$ and 
$(\la_\dashv)^\circ=(\la^\circ)_\vdash$, it follows that 

(d) $(\boldkey a_\vdash)'=(\boldkey a')_\dashv$ and 
$(\boldkey a_\dashv)'=(\boldkey a')_\vdash$.
\enddefinition

Let $\fS_r=\fS_{\{1,\cdots,r\}}$ be the symmetric group
on $r$ letters as in the introduction. 
Each element $w\in \frak S_r$ can be expressed as a 
product of $s_{i_1}\cdots
s_{i_k}$, where $s_i=(i,i+1)$ are basic transpositions. If
 $k$ is minimal, such an expression is called a reduced expression
of $w$. The number $k$ is defined to be the length $l(w)$ of $w$. It is
independent of the reduced expression of $w$ and
$l(w)=\# \{ (i, j)\mid i<j, (i) w> (j) w\}.$
For  $\boldkey a=[a_i]\in \La[m, r]$ with $\la=\Theta(\boldkey a)$,  let 
$$
\frak S_{\boldkey a}=\fS_\la=\frak S_{\{1,\cdots,a_1\}}\times 
\frak S_{\{a_1+1, \cdots, a_2\}}
\times \cdots \times  \frak S_{\{a_{m-1}+1, \cdots, a_m\}}.
$$
be  the Young subgroup of $\frak S_r$ corresponding to $\bka$,
and $\sD_\bka=\sD_\la$ the set of distinguished representatives of right
$\frak S_\la$-cosets. 
Similarly, for $\bka,\bkb\in\La[m,r]$,
  $\sD_{\bka,\bkb}=\sD_\bka\cap \sD_\bkb^{-1}$ denotes the set of distinguished
representatives of  $\frak S_\bka$-$\frak S_\bkb$ double cosets.

For positive integers  $i, j$, let $s_{i,i}=1$ and
$$
s_{i, j}=\cases s_{i-1}\cdots s_j, &\text{if $i>j$,}\cr
                 s_{i}\cdots s_{j-1}, &\text{if $i<j$.} \cr
\endcases$$

\proclaim{(1.3) Lemma} (a) We have $\frak S_r=\cup_{i=1}^r s_{i, r} \frak S_{r-1}$,
where $s_{i, r}$, $1\le i\le r$, are 
 distinguished $\fS_{r-1}$-coset representatives in $\fS_r$.

(b) Let $\bka\in \La[m, r]$ and $d\in \sD_{\bka, \bka'}$.
Then  
$d=s_{a_j, r} d_1$ for some $ j$ with $a_{j-1} < a_j$ and   
$d_1\in \sD_{\bkb, (\bka_\vdash)'}$, where  
$\bkb=[a_0, \cdots, a_{j-1}, a_j-1, a_{j+1}-1, 
\cdots, a_m-1]$. 
\endproclaim

\demo{Proof} The statement  (a) is well-known. Since $d\in \sD_{\bka'}^{-1}$, we have
$d_1\in \sD_{\bka'}^{-1}\cap \fS_{r-1}$. Therefore, 
$d_1\in \sD_{(\bka')_\dashv}^{-1}=\sD_{(\bka_\vdash)'}^{-1}$ by (1.2)(d).
On the other hand, $d\in \sD_{\bka}$ implies that $i=a_j$ for some 
$1\le j\le m$. Take the minimal $j$ with $a_j=i$. Then 
 $a_{j-1}\neq a_j$, and $\bkb\in \La[m, r-1]$.
If $d_1\not \in \sD_{\bkb}$, then there is a $s_k\in \fS_\bkb$ such that
$d_1=s_k d_2$ with $l(d_1)=l(d_2)+1$, and 
$d=s_{a_j, r} s_k d_2$. If $k\ge a_j$, then 
$d=s_{k+1} s_{a_j, r}  d_2$, contrary to  $d \in \sD_\bka$.
If $k<a_j$, then $k\le a_j-2$ since $s_k\in \frak S_{\bkb}$.
Thus,  $d=s_{k} s_{a_j, r}  d_2$, a contradiction again. Therefore, 
 $d_1\in \sD_{\bkb}$, and $d_1\in \sD_{\bkb, (\bka_\vdash)'}$. 
\qed\enddemo

Let $w_{i, j}  =s_{i+1, 1} s_{i+2, 2}\cdots s_{i+j, j}$. Then $w_{i,j}$
is the following permutation 
$$w_{i,j}=\left(\smallmatrix 1&\cdots&i&i+1&\cdots&i+j\cr
j+1&\cdots&j+i&1&\cdots&j\cr\endsmallmatrix\right).
\tag1.4$$
Let $k$ be a non-negative integer. Define the $k$-shifted elements
$s_{i, j}^{(k)},  w_{i, j}^{(k)}$ by setting 
$s_{i, j}^{(k)}=s_{i+k, j+k}$ and $w_{i, j}^{(k)} =s_{i+1, 1}^{(k)} 
s_{i+2, 2}^{(k)} \cdots s_{i+j, j}^{(k)}$, and
 $w_{i,j}=w_{i,j}^{(k)}=1$, if $i=0$ or $j=0$.
Note that $w_{i, j}^{(k)}$ is a permutation on
$\{k\!+\!1,\cdots,k\!+\!i\!+\!j\}$, and explicitly,
$$w_{i,j}^{(k)}=\left(\smallmatrix k+1&\cdots&k+i&k+i+1&\cdots&k+i+j\cr
k+j+1&\cdots&k+j+i&k+1&\cdots&k+j\cr\endsmallmatrix\right).\tag1.5$$
Obviously, we have 
$(s_{i, j}^{(k)})^{-1}=s_{j, i}^{(k)}\text{ and } 
(w_{i, j}^{(k)})^{-1}=w_{j, i}^{(k)}$

For $\boldkey a=[a_i]\in \La[m, r]$ with $\la=\Theta(\bka)$, let  
 $w_{\boldkey a}\in \frak S_r$  be defined by 
$$
(a_{i-1}+l) w_{\boldkey a}=r-a_i+l\text{ for all }i
\text{ with } a_{i-1}<a_i, 1\le l\le a_{i}-a_{i-1}.\tag1.6
$$ 
In particular, we have $(a_i) w_{\bold a} =r-a_{i-1}$ if $a_{i-1}<a_i$.
For example, for
$\bka=[0,i,i+j]$, $w_\bka=w_{i,j}$, and for $\bkb=[0,2,5,9]$,
$$w_{\bkb}=\left(\smallmatrix 1\,\,2&3\,\,4\,\,5&6\,\,7\,\,8\,\,9\cr
                              \underline{8\,\,9}&\underline{5\,\,6\,\,7}
&\underline{1\,\,2\,\,3\,\,4}\cr
\endsmallmatrix\right).$$

\proclaim{(1.7) Lemma} For $\boldkey a=[a_i]\in \La[m, r]$ with $\la=\Theta(\bka)$, 
$$\aligned
w_{\boldkey a} &= w_{a_{m-1}, \la_m}^{(0)} 
w_{a_{m-2}, \la_{m-1}}^{(\la_m)} w_{a_{m-3}, \la_{m-2}}^{(\la_m+\la_{m-1})}
\cdots w_{a_{1}, \la_2}^{(\la_m+\cdots+\la_3)}\cr
&= w_{a_{m-1}, a_m-a_{m-1}}^{(0)} w_{a_{m-2}, a_{m-1}-a_{m-2}}^{(a_m-a_{m-1})}
\cdots w_{a_{1}, a_2-a_{1}}^{(a_m-a_2)}.\cr
\endaligned$$ 
\endproclaim
\demo{Proof} The result follows immediately from (1.4) and (1.5). 
\qed\enddemo

We list some properties for the elements $w_\bka$. First,
$w_\bka$ is distinguished and turns $\fS_\bka$ into $\fS_{\bka'}$.

\proclaim{(1.8) Lemma}  For any  $\boldkey a=[a_i]\in \La[m, r]$, we have

(a) $w_{\boldkey a}^{-1} =w_{ \boldkey a'}$ and $w_\bka\in\sD_\bka$. Hence
$w_{\boldkey a}\in \sD_{\boldkey a,\bka'}$.

(b) $w_{\boldkey a}^{-1} \frak S_{\boldkey a} w_{\boldkey a}=\frak S_{\boldkey a'}$.
In particular, $w_{\boldkey a}^{-1}s_j w_{\boldkey a}=s_{(j)w_{\boldkey a}}$
for all $j\ge 1$ with $j\neq a_i$.
\endproclaim

\demo{Proof} The first assertion in (a) follows from definition.
Consider the root system $\Phi$ of type $A_{r-1}$ and its 
subsystem $\Phi_\bka$ whose Coxeter graph  is obtained by removing all $a_i$-th
vertices from that of $\Phi$.
For any $s_i\in \frak S_r$, there is   a simple root
$\alpha_i=e_i-e_{i+1}$ with $s_i=s_{\alpha_i}$ and
$w(\alpha_i)=e_{(i) w^{-1}}-e_{(i+1) w^{-1}}$ (see, e.g.,  \cite{Hum}). 
From  (1.6),  we have $(a_j) w_{\boldkey a}=r-a_{j-1}$ for $a_{j-1}<a_j$, 
and $(i+1) w_{\boldkey a}=(i)w_{\boldkey a}+1$
if  $i\neq a_j$ for all  $1\le j\le m$. Thus, $w_{\boldkey a}^{-1} (\alpha_i)$
is a simple root.
Therefore, $w_{\bka}^{-1}$ stabilizes the positive root system of $\Phi_\bka$,
and consequently, $ w_{\boldkey a} \in \sD_{\bka}$  
and  hence $w_{\boldkey a} \in \sD_{\bka, \bka'}$, proving (a). 
Because 
 $w_{\boldkey a}^{-1}  s_i w_{\boldkey a}=s_{(i)w_{\boldkey a}}
\in \frak S_{\boldkey a'}$ is  a basic transposition, 
 $w_{\boldkey a}^{-1}   \frak S_{ \boldkey a} w_{\boldkey a}
\subseteq \frak S_{\boldkey a'}$.  Therefore, by (a),
$w_{\boldkey a}^{-1}   \frak S_{ \boldkey a} w_{\boldkey a}=\frak S_{ \boldkey a'}$, proving (b).
\qed\enddemo

We now look at the relation between $w_\bka$ and $w_{\bka_\dashv}$.
 For $\boldkey a=[a_i]\in \La[m, r]$ with minimal index $k$ such that 
$a_k=r$, we have (see (1.2(c)) 
$\bka_\dashv=[0, a_1, \cdots, a_{k-1},  r-1, \cdots, r-1]\in \La[m, r-1]$.
Let  $\boldkey a_i\in \La[m, r]$, $i=1,\cdots,m$, be defined by 
 $$
\bka_i=\bka_\dashv+{\boldkey 1}_i,
\text{ where }{\boldkey 1}_i=[0,\undersetbrace{i-1}\to{0,\cdots,0},1,\cdots,1]\in
\La[m,1]. 
\tag1.9
$$
Then $\boldkey a_i\in \La[m, r]$ with 
 $\boldkey a_1\rec \boldkey a_2 \rec \cdots\rec \bka_m
$ and  $\bka_k=\bka$.

\proclaim{(1.10) Lemma} Write $\bka_\dashv=[b_0, b_1, \cdots, b_m]\in \La[m, r-1]$.
Then, for any $i$, $1\le i\le m$, 
$w_{\boldkey a_i} = s_{b_i+1, r} w_{\boldkey a_\dashv} s_{r, r-b_{i-1}}$
with $l(w_{\boldkey a_i})=l(s_{b_i+1, r} )+l(w_{\boldkey a_\dashv})+l(s_{r, r-b_{i-1}})$.
\endproclaim

\demo{Proof}   We prove $
(l) w_{\boldkey a_i} =
 (l) s_{b_i+1, r} w_{\boldkey a_\dashv} s_{r, r-b_{i-1}}
$ 
for every  $l$ with $1\le l\le r$. 

Assume that $l\le b_{i}$. Then there is a $j\le i$ with $b_{j-1}  < l\le b_j$. 
Write  $l=b_{j-1}+l'\le b_{j}$. 
Then  $(l) s_{b_i+1, r} w_{\boldkey a_\dashv} s_{r, r-b_{i-1}} =(l)
 w_{\boldkey a_\dashv} s_{r, r-b_{i-1}}
=((r-1)-b_j+l') s_{r, r-b_{i-1}}$ using (1.6).
Because  $r-1-b_j+l' \ge r-b_{i-1}$ for $j<i$ and $r-1+l'-b_i< r-b_{i-1}$,
 we have  
$$
(l) s_{b_i+1, r} w_{\boldkey a_\dashv} s_{r, r-b_{i-1}}=\cases r-b_j+l'  & \text{ if $j<i$}\cr
                                   (r-1)-b_j+l' & \text { if $j=i$.}\cr\endcases 
$$
Thus, $ (l) s_{b_i+1, r} w_{\boldkey a_\dashv}
 s_{r, r-b_{i-1}}=(l)w_{\boldkey a_i}$. 
If $l=b_i+1$, then $(l) s_{b_i+1, r} w_{\boldkey a_\dashv} s_{r, r-b_{i-1}}=r-b_{i-1}
=(l)w_{\boldkey a_i}$.
 Assume  $l>b_i+1$.
Then there is a $j$ with $j\ge i$ and $b_j< l-1\le b_{j+1}$. 
Write $l-1=b_j+l'$.
Because $r-1-b_{j+1}+l'<  r-b_{i-1}$, we have 
$(l) s_{b_i+1, r} w_{\boldkey a_\dashv} s_{r, r-b_{i-1}} =
 (l-1) w_{\boldkey a_\dashv} s_{r, r-b_{i-1}}
  =r-1-b_{j+1}+l'$.  In this case, $l=(b_j+1)+l'$ and 
$(l)w_{\boldkey a_i}=r-(b_{j+1}+1)+l'=
(l) s_{b_i+1, r} w_{\boldkey a_\dashv} s_{r, r-b_{i-1}} $. So 
$w_{\boldkey a_i}=s_{b_i+1, r}w_{\boldkey a_\dashv} s_{r, r-b_{i-1}}$.
The length formula is obviously.
\qed \enddemo

The elements $w_\bka\in\fS_r$ and $w_{\bka_\vdash}\in\fS_{r-1}$ are related as follows.

\proclaim{(1.11) Corollary} Let $\boldkey a=[a_i]
\in \La[m, r]$ and
suppose that $k$ is the minimal index with $a_k\neq 0$. Then,  

(a)  $w_{\boldkey a} = s_{a_k, r} w_{\boldkey a_\vdash }$ with $l(w_{\boldkey a})
=l(w_{\boldkey a_\vdash })+ l(s_{a_k, r})$. 

(b)  $w_{\boldkey a'}=w_{(\boldkey a_\vdash)'} s_{r, a_{k}}$ with 
$l(w_{\boldkey a'})=l( w_{(\boldkey a_\vdash)'})+l(s_{r, a_{k}})$.
\endproclaim

\demo{Proof}  It follows immediately from (1.2)(d) and 
(1.10), or from (1.8) and (1.3).\qed\enddemo

\subhead 2. Some zero divisors\endsubhead
Let $W=W_m^r$ be the group defined in the introduction. Then $W$ is the  group
with $r$ generators 
$\{s_0, s_1, s_2, \cdots, s_{r-1}\}=S$ and
relations:
$$\aligned
&  s_0^m =s_i^2=1, \text{ for all $1\le i\le r-1$,}\cr
&  s_is_{i+1}s_i=s_{i+1}s_is_{i+1}, \text{ for $1\le i\le r-2$,}\cr
&  s_0s_1s_0s_1=s_1s_0s_1s_0 \cr
&  s_is_j=s_js_i, \text{ if $0\le i\le j-2\le r-3$}.\cr 
\endaligned
$$
Let $t_1=s_0$ and $t_i=s_{i-1}t_{i-1}s_{i-1}$,  $2\le i\le r$.
Then $t_it_j=t_jt_i$, $t_i^m=1$ for  $1\le i, j\le r$
and $\{t_i\mid  1\le i\le r\}$ generates the bottom $C$.
We will identify the subgroup of $W$  generated by 
$\{s_i\mid 1\le i\le r-1\}$ with $\frak S_r$. 
Note that the group $W_m^r$ is the symmetric 
(resp. hyperoctahedral) group if $m=1$ (resp. $m=2$).

A deformation of the group algebra of $W$ has been  given recently by Ariki and 
Koike \cite{AK}. Let $R$ be a commutative ring with  $1$ and  
$q, q^{-1}, u_1, \ldots, u_m\in R$. The Ariki-Koike algebra $\bH=\bH_R=
\bH_m^r$ associated
to the group $W$ is an  associative algebra over $R$ 
with generators $T_i:=T_{s_i}$,  $0\le i\le r-1$, subject to the  relations:  
$$\cases 
          (1) &T_0T_1T_0T_1 = T_1T_0T_1T_0, \cr  
           (2) & T_iT_{i+1}T_{i}=T_{i+1}T_iT_{i+1}, \hskip.5cm  
         \text{ for $1\le i\le r-2$}\cr 
          (3)& T_i T_j=T_jT_i,\hskip.5cm  \text{ if  $|i-j|\ge 2$}\cr 
          (4)& (T_i-q)(T_i+1)=0,\hskip.5cm  \text{ if  $i\neq 0$}\cr 
          (5)& (T_0-u_1)\cdots (T_0-u_m)=0.\cr
\endcases \tag 2.1
$$
 Let $L_1=T_0$ and  $L_i=q^{-1} T_{i-1}L_{i-1}T_{i-1}$, $2\le i\le r$. Then
$L_iL_j=L_jL_i$, $1\le i, j\le r$. The elements $L_i$, $1\le i\le r$,
generate an abelian subalgebra of $\bH$. Since $L_i^m\not=1$ in general,
this subalgebra contains a proper submodule of rank $|C|$. So it cannot
 serve as a ``bottom'' of $\bH$.
However, $\bH$ is $R$-free of rank $|W|$ with  basis  \cite{AK}:
$$
\{L_1^{c_1}\cdots L_r^{c_r} T_w\mid w\in \frak S_r, \text{ and } 0\le c_i\le 
m-1, \forall  i\}.\tag 2.2 
$$
Let $\sH$ be the subalgebra generated by $T_i$ for all $i\ge1$. 
For a Young subgroup $\frak S_{\bka}$, let $\sH(\frak S_{\bka})$ be the
corresponding subalgebra. By (2.2), we have  
$\bH=\oplus_{\bold c}L^{\bold c}\sH$, where $\bold c=(c_1,\cdots,c_r)$
and $L^{\bold c}= L_1^{c_1}\cdots L_r^{c_r}$. Thus, for every such $\bold c$,
we have a projection map
$$pr_{\bold c}:\bH\to L^{\bold c}\sH
.\tag2.3$$

The isomorphism from $\bH$ to the algebra $\bH^\op$ opposite to
$\bH$ induces an
anti-automorphism
$$\iota: \bH\rightarrow \bH \text{ such that }  \iota(T_i)=T_i. \tag 2.4
$$
Clearly, $\iota(L_i)=L_i$. We need the following commutator relations.

\proclaim{(2.5) Proposition}
Let $\bH$ be the Ariki-Koike  algebra over a commutative 
ring  $R$.  
Then 

(a) $T_i$ commutes with $L_j$ if  $j\neq i, i+1$.

(b) $T_i$ commutes with $L_iL_{i+1}$ and $L_i+L_{i+1}$.  

(c) $T_i$ commutes with $\prod_{j=1}^k (L_j-x)$ for all $x\in R$ and $i\neq k$.

(d)  $ L_i^k=q^{-1}T_{i-1} L_{i-1}^k T_{i-1} +(1-q^{-1})\sum_{c=1}^{k-1}
L_i^cL_{i-1}^{k-c} T_{i-1}$ if $1<i\le r$ and $ k\ge 1 $.
\endproclaim

\demo{Proof} See \cite{AK, (3.3)} for statements (a-c), and
\cite{MM, (3.6)} for (d) .
\qed\enddemo

Following Graham and Lehrer \cite{GL, (5.4) } (or \cite{DJM, (3.1)}), we
define for $\bka=[a_i]\in\La[m,r]$,   
$$\cases
 \pi_{\boldkey a}=\pi_{a_1}(u_2)   \cdots \pi_{a_{m-1}}(u_m),\qquad
\tpi_{\boldkey a}=\pi_{a_1} (u_{m-1})   \cdots \pi_{a_{m-1}} (u_1),& \cr
 \text{where } \pi_0(x)=1 \text { and } 
\pi_{a} (x)= \prod_{j=1}^{a} (L_j -x),\forall a>0, x\in R.&
\cr\endcases \tag2.6
$$

Note that $\pi_\bka=\pi_{\bka_\dashv}$ if $a_{m-1}\neq r$ and
$\tpi_{\bka'}=\tpi_{(\bka_{\vdash})'}$ if $a_1\neq 0$.
Also, if, for  $x_1,\cdots,x_{m-1}\in R$, we put
$\pi(\bka;x_1,\cdots,x_{m-1})=\pi_{a_1}(x_1)   \cdots \pi_{a_{m-1}}(x_{m-1})$,
then $\pi_\bka=\pi(\bka;u_2,\cdots,u_{m})$ and
$\tpi_\bka=\pi(\bka;u_{m-1},\cdots,u_{1})$.

\proclaim{(2.7) Corollary} (a)  For $\boldkey a\in \La[m, r]$, 
$\pi_{\boldkey a}$ and 
$\tpi_{ \boldkey a}$ commute with any element in $\sH(\frak S_{\boldkey a})$.
In particular, for any $x\in R$, $\pi_r(x)$ is in the centre
of $\bH$.

(b) Assume $a_{j-1}<i\le a_j$ for some $j$. Then   
$$
\pi_{\boldkey a} T_{i, r} =T_{i, a_j} (L_{a_j}-u_{j+1}) T_{a_j, a_{j+1}}
\cdots (L_{a_{m-1}}-u_m) T_{a_{m-1}, a_m}  
 \pi_{\boldkey b},
$$
where $\boldkey b=[0, \cdots, 0, a_k, \cdots, a_{j-1}, a_j-1, \cdots,
a_m-1]\in \La[m, r-1]$ and $T_{i,j}=T_{s_{i,j}}$. 
\endproclaim 

\demo{Proof} The statement (a) follows from (2.5)(c). 
Noting that $T_{i, j}=T_i\cdots T_{j-1}$ for $i<j$, we obtain (b)
immediately from (2.5)(a)-(c).
\enddemo

We now prove the following useful  result.

\proclaim{(2.8) Lemma} For  $\boldkey a, \boldkey b\in \La[m, r]$,
we have
$\pi_{\boldkey a} \sH \tpi_{\boldkey b'}  =0$ 
and $\tpi_{\boldkey a} \sH
\pi_{\boldkey b'}=0$
unless $\boldkey a\lecq \boldkey b$.
\endproclaim

\demo{Proof} Using the anti-automorphism $\iota$ of $\bH$ in (2.4)
and the fact $\boldkey b'\lecq \boldkey a'$ if and only if $\boldkey a\lecq 
\boldkey b$,  we see that both assertions in (2.8) are equivalent. 
Therefore, we only need to  prove $\pi_{\boldkey a} \sH \tpi_{\boldkey b'}  =0$
 unless  $\boldkey a\lecq \boldkey b$. 

We apply induction on $r$.
Let $r=1$. If $\boldkey a\not\lecq \boldkey b$, then there is an
$i$ with $a_{i}>b_{i}$. Since $b_i, a_i\le 1$ for all $i$, we have 
$b_i=0$ and $a_i=1$. By (2.6), $\prod_{k=i}^{m-1}  (L_1-u_{k+1})
=\prod_{k=i+1}^m (L_1-u_k)$ 
( resp. $\prod_{k=1}^i  (L_1-u_{k})$ )
is a factor of $\pi_{\boldkey a}$ (resp. $\tpi_{\boldkey b'}$).  
Therefore, $\pi_{\boldkey a} \sH \tpi_{\boldkey b'}  =0$
by (2.1)(5), and (2.8) is true for $r=1$. 
Assume $\pi_{\boldkey a} \sH\tpi_{\boldkey b'}  =0$
for all $\bka, \bkb \in \La[m, r-1]$ with $\bka\not\lecq \bkb$.

Let  $i$ and $j$ be the minimal indices  with 
$b_i\neq 0$ and $a_j\neq 0$,  respectively.  Because $b_k=0$ for
all $k<i$, 
$\pi_{b_m-b_k}(u_k) =\pi_r(u_k)$, which are in the centre of $\bH$ (see
(2.7)(a)).
 Therefore,  for any $w\in \frak S_r$,  $\pi_{\boldkey a} 
T_w \tpi_{\boldkey b'} $ contains a factor
$\prod_{k=j+1}^m (L_1-u_k) \prod_{k=1}^{i-1} (L_1-u_k)$. 
By (2.1)(5),  $\pi_{\boldkey a} \sH \tpi_{\boldkey b'}=0$  
 unless $i\le j$.  On the other hand, take $w\in \frak S_r$ 
with $\pi_{\boldkey a} T_w \tpi_{\boldkey b'} \neq0$.
 Write $w=dy$ with $y\in \frak S_{r-1}$ and  $d=s_{k, r}$ 
for some $1\le k\le r$ (see (1.3)). 
By (2.7)(b),  $\pi_{\boldkey a} T_d = h\pi_{\boldkey a_\vdash}$ for some 
$h\in \bH$.
By (2.6),  $\tpi_{ \boldkey b'} = \tpi_{ (\bkb_\vdash)'} 
(L_r-u_{i-1})\cdots (L_r-u_1)$. 
So, $\pi_{\boldkey a} T_w \tpi_{\boldkey b'} \neq0$ implies
 $ \pi_{\boldkey a_\vdash } T_y  \tpi_{(\boldkey b_\vdash)'} \neq 0$.
Now, by induction, 
  $\boldkey a_\vdash \lecq \boldkey b_\vdash$, 
which implies $\boldkey a\lecq \boldkey b$ since  $i\le j$.
\qed\enddemo

\subhead 3. Idempotents\endsubhead In this section,  
idempotents $e_{\boldkey a}$ for all $\boldkey a\in \La[m, r]$
will be constructed under a certain condition.

\proclaim{(3.1) Proposition}   For any $\boldkey a\in \La[m, r]$, 
let $v_{\boldkey a}=\pi_{\boldkey a} T_{w_{\boldkey a}}\tpi_{\boldkey a'}$.
Then we have:

(a) $\pi_{\boldkey a} \sH\tpi_{\boldkey a'}=v_\bka\sH(\fS_{\bka'})=\sH(\fS_{\bka})v_{\bka}$.
Moreover, 
$v_{\boldkey a} T_i= T_{(i) w_{\boldkey a}^{-1} }v_{\boldkey a}$ for 
any $s_i\in \frak S_{\boldkey a'}$.

(b) $v_{\boldkey a} L_i \in v_{\boldkey a} \sH(\frak S_i)\cap 
v_{\boldkey a} \sH(\frak S_{\bka'})$ and 
$ L_iv_{\boldkey a}\in  \sH(\frak S_i) v_{\bka}\cap \sH(\frak S_{\bka})
v_{\boldkey a}$ 
for every $i=1,\cdots,r$.
In particular,  $v_{\boldkey a} L_i =u_j v_{\boldkey a}$ if 
$i=r-a_j+1$ for some $j$ with $a_{j-1}<a_j$.

(c) $\pi_{\boldkey a}\bH \tpi_{\boldkey a'}=\pi_{\boldkey a}\sH\tpi_{\boldkey a'}$.

(d) $v_{\boldkey a} \bH=v_{\boldkey a} \sH$.
\endproclaim

\demo{Proof} Let $d\in \sD_{\bka, \bka'}$. We prove  
  $\pi_{\bka} T_d \tpi_{\bka'}= 0$ 
unless $d=w_{\bka}$. Obviously, 
this result is true for $r=1$, and assume that $r>1$ and that 
the result is true
for $r-1$.
By (1.3),  $d=s_{i, r} d_1$ with 
$d_1\in \frak S_{r-1}$ and $i\geq a_k$, where
$k$ is the minimal index with $a_k\neq 0$. 
If $i>a_k$, then    there is an index  $j$ with 
$a_{j-1}<i\le a_{j}$, $j>k$. By (2.7)(b), 
$$
\pi_{\boldkey a} T_d =
T_{i, a_j} (L_{a_j}-u_{j+1}) T_{a_j, a_{j+1}}
(L_{a_{j+1}}-u_{j+2})\cdots (L_{a_{m-1}}-u_m) T_{a_{m-1}, a_m}  
 \pi_{\boldkey b} T_{d_1},
$$
where $\boldkey b=[0, \cdots, 0, a_k, \cdots, a_{j-1}, a_j-1, \cdots,
a_m-1] \in \La[m, r-1]$ and  $\boldkey b\not\lecq \boldkey
a_\vdash$. By (2.6),  $\tpi_{\boldkey a'}=\tpi_{  (\boldkey a_\vdash) '} (L_r-u_{k-1})\cdots
(L_r-u_1)$. We have $\pi_\boldkey b T_{d_1} \tpi_{  (\boldkey a_\vdash)'}=0$
by (2.8). So    $\pi_{\boldkey a} T_d
\tpi_{\boldkey a'}= 0$, a contradiction, proving  
 $i=a_k$. From (1.3)(b),  we have $d_1\in \sD_{\bka_\vdash, (\bka_\vdash)'}$. Because 
$$ \aligned 
\pi_{\bka} T_d \tpi_{\bka'} = & (L_{a_k}-u_{j+1}) T_{a_k, a_{k+1}}
(L_{a_{k+1}}-u_{k+2})\cdots (L_{a_{m-1}}-u_m) T_{a_{m-1}, a_m} \cr
 &  \pi_{\boldkey a_\vdash} T_{d_1} \tpi_{  (\boldkey a_\vdash) '} 
(L_r-u_{k-1})\cdots (L_r-u_1), \cr\endaligned
$$ 
we have $\pi_{\boldkey a_\vdash} T_{d_1} \tpi_{(\boldkey a_\vdash) '}\neq 0$. 
By induction,  $d_1=w_{\bka_\vdash}$, and therefore,  $d=w_{\bka}$ by (1.11).
Now, the first assertion in (a)  follows immediately from the  
$\frak S_{\bka}$-$\frak S_{\bka'}$
double coset decomposition   of $\frak S_r$.

For notational simplity, we write 
$$h_{\bka}=T_{w_\bka}.$$   By (1.8), we have  
$h_{\boldkey a} T_i= T_{(i) w_{\boldkey a}^{-1} } h_{\boldkey a} $ for every $s_i
\in \frak S_{\boldkey a'}$. Therefore, 
$v_{\boldkey a} T_i=\pi_{\boldkey a} h_{\boldkey a} T_i \tpi_{{\boldkey a}'} 
=\pi_{\boldkey a} T_{(i) w_{\boldkey a}^{-1} } h_{\boldkey a} 
 \tpi_{{\boldkey a}'} = T_{(i) w_{\boldkey a}^{-1} }v_{\boldkey a}$, proving 
the second assertion in (a).

To see (b), we first treat the case  $i=a_m-a_j+1$ with $j$ minimal. Then  $a_j>a_{j-1}$, and 
   $\boldkey b=[a_0, a_1, \cdots, a_{j-1}, a_j-1, 
a_{j+1}, \cdots, a_m]\in \La[m, r]$.   
Obviously, $\boldkey b\not\recq \boldkey a$. By (2.8),    
$\pi_{\boldkey a}  h_{\boldkey a} \tpi_{\boldkey b'} =0$. 
By (2.6), we have $\tpi_{\boldkey a'} L_i =\tpi_{\boldkey b'}+u_j\tpi_{\boldkey a'}$,
where 
$\boldkey b'=[0, a_m-a_{m-1}, \cdots, a_{m}-a_{j-1}, a_{m}-(a_j-1), \cdots, a_m-a_1, a_m].
$ 
Thus, 
$v_{\boldkey a} L_{i}=\pi_{\boldkey a}  h_{\boldkey a} \tpi_{\boldkey b'}+
u_{j} v_{\boldkey a} =u_{j} v_{\boldkey a}$.
Now, assume $i\neq a_m-a_j+1$, $1\le j\le m$. Then there is a $k$ with 
$a_m-a_{k}+1< i\le a_m -a_{k-1}$. Write $i=(a_m-a_k+1)+l$. Then 
$L_i=q^{-l} T_{i, a_m-a_k+1} L_{a_m-a_k+1}  T_{ a_m-a_k+1, i}$.
By (a), $v_{\boldkey a} L_i=u_k q^{-l} v_\bka T_{i, a_m-a_k+1}  T_{ a_m-a_k+1, i}\in 
v_{\boldkey a} \sH(\frak S_i)\cap  v_{\boldkey a} 
 \sH(\frak S_{\boldkey a'})$. By a symmetric argument, 
one can prove $L_iv_{\boldkey a}\in  \sH(\frak S_i)
 v_{\bka}\cap \sH(\frak S_{\bka})
v_{\boldkey a}$. 

By (b) and (a), $L_i v_{\bka}\sH(\frak S_{\bka'})\subseteq 
\sH(\frak S_{\bka})v_{\bka}\sH(\frak S_{\bka'})
= v_{\bka}\sH(\frak S_{\bka'})
$ for every $i=1,\cdots, r$. 
Thus, (c) follows immediately from (2.2). 

For arbitrary $v_{\boldkey a}T_w \in v_{\boldkey a}\sH$
with $w\in \frak S_r$, write $w=dy$ with $y\in \frak S_{\{2, \cdots, r\}}$ and
$d\in \frak S_r/ \frak S_{\{2, \cdots, r\}}$. Then $d=s_{i, 1}$
with $1\le i\le r$. By (b), 
$v_{\boldkey a} T_w T_0=v_{\boldkey a} T_{i, 1} T_0 T_y\in v_{\boldkey a}\sH$. 
Thus,  $v_{\boldkey a}  \sH$ is stable under the right action of
$T_0$, and  
$v_{\boldkey a} \bH=v_{\boldkey a} \sH$.
\qed\enddemo

Recall from (2.7) that $T_{i,j}=T_{s_{i,j}}$ for $i\neq0\neq j$.
 For $\bka\in \La[m, 0]$, let  $v_{\bka}=1$.

\proclaim{(3.2) Lemma}  
Let $x_1,\cdots,x_{m-1}$ be elements in the commutative ring $R$
and $\bka\in \La[m, r]$ with $a_1\neq 0$.
Then there is an integer $c$ depending on $\bka$ such that 
 $$
T_{a_m, a_{m-1}} (L_{a_{m-1}}-x_{m-1}) \cdots T_{a_2, a_1}  (L_{a_1}-x_1)  
  -q^{c} L_r^{m-1} T_{a_1, a_{m}}^{-1} 
$$
is in the free $R$-submodule spanned by 
$\{  L_1^{c_1}\cdots L_r^{c_r}  T_w \mid w\in \frak S_r, c_j< m-1, \forall j\}$.
\endproclaim

\demo{Proof} Let $U_i$ ($1\le i\le m-1$) be the free $R$-submodule spanned by 
$$\{ L_1^{c_1}\cdots L_r^{c_r} T_w \mid w\in \frak S_{\{a_i,  a_i+1, \cdots, r\}}, 
c_j<m-i, \forall j\}.
$$
We claim that, for any $i=1,\cdots,m-1$,  
$$ T_{a_{m}, a_{m-1}}(L_{a_{m-1}}-x_{m-1})\cdots 
T_{a_{i+1}, a_i} (L_{a_i}-x_{i})
  -q^{c}  L_r^{m-i} T_{a_i, a_{m}}^{-1}\in U_i.\tag 3.3 
$$ 
Apply induction on $i$.
The result for $i=m-1$ is true
 since  
$T_{a_{m}, a_{m-1}} (L_{a_{m-1}}-x_{m-1}) =q^{a_{m}-a_{m-1}} L_{a_{m}}  
T_{a_{m-1}, a_{m}}^{-1} 
-x_{m-1} T_{ a_{m}, a_{m-1}}$. 
For  $i<m-1$, we have,
by induction,  
$$\aligned &
T_{a_m, a_{m-1}}(L_{a_{m-1}}-x_{m-1})\cdots T_{a_{i+1}, a_i} (L_{a_i}-x_i) \cr = &
(q^{c} L_r^{m-i-1} T_{a_{i+1}, a_m}^{-1}  +h) T_{a_{i+1}, a_i} (L_{a_i}-x_{i})\cr
=& q^{c} L_r^{m-i-1} T_{a_{i+1}, a_m}^{-1} T_{a_{i+1}, a_i} (L_{a_i}-x_{i})
+ h T_{a_{i+1}, a_i} (L_{a_i}-x_{i})\cr\endaligned
$$
for some $h\in U_{i+1}$. We may assume $h=L_1^{c_1}\cdots L_r^{c_r} T_w $ 
with  $ w\in \frak S_{\{a_{i+1},  \cdots, r\}}$ and  $c_j<m-(i+1)$, $\forall j$,
without loss of generality. Write $w=s_{k, a_{i+1}} y$ with 
$y\in \frak S_{\{a_{i+1}+1, \cdots, r\}}$ (cf. (1.3)). Then
$T_w T_{a_{i+1}, a_i}(L_{a_i}-x_{i+1})=q^{k-a_i}L_k T_{a_i, k}^{-1}T_y 
-x_{i+1} T_{k, a_i}T_y$.
Thus, we have  $h T_{a_{i+1}, a_i} (L_{a_i}-x_{i+1})\in U_i$.
Noting  $qT_i^{-1}=T_i-(q-1)$ and 
$L_r^{m-i-1}T_{r, a_{i+1}} T_{ a_{i+1}, a_i} L_{a_i}
=q^{r-a_i} L_r^{m-i}T_{ a_i, r}^{-1}$, we may write
$$
L_r^{m-i-1} T_{a_{i+1}, a_m}^{-1} T_{a_{i+1}, a_i} (L_{a_i}-x_{i+1})
=q^{c'} L_r^{m-i} T_{a_i, r}^{-1}+h'
$$  
for some  $h'\in U_i$ and  $c'\in \Bbb Z$.
Therefore, (3.3) holds for $i$, proving the claim.
Now, the required result follows from the case  $i=1$.
\qed\enddemo

\proclaim{(3.4) Theorem} Let $\bH$ be the Ariki-Koike  algebra 
over $R$. For any $\boldkey a\in \La[m, r]$, 
$v_{\boldkey a} \bH$ is a free
$R$-module  with basis $\{v_{\boldkey a} T_w\mid w\in \frak S_r\}$.
\endproclaim

\demo{Proof} For $h=(m-1, \cdots, m-1)$, let $pr_h $ be the projection defined in (2.3).
First, we prove 
$$
pr_h (v_{\boldkey a})=q^c \prod_{j=1}^{r}L_j^{m-1} h_{\boldkey a'}^{-1}
\tag 3.5
$$
for some  integer $c\in \Bbb Z$ depending on $\bka$. We prove (3.5)
by induction on $r$. Obviously, the result is true for $r=1$. Assume
now that $r>1$ and that (3.5) holds for $r-1$.
For $\bka\in \La[m, r]$ with minimal index $k$, $a_k= r$,  recall
$$\aligned & \boldkey a'=[0, 0, \cdots,0, r-a_{k-1}, \cdots,
r-a_1, r],\cr
           & \bka_\dashv=[0, a_1, \cdots, a_{k-1}, r-1, \cdots, r-1],\cr
           & (\bka_\dashv)'=[0, \cdots, 0, r-1-a_{k-1}, \cdots, r-1-a_1, r-1],
\text{ and }\cr
& \pi_{\bka}=\pi_{\bka_\dashv}(L_r-u_{k+1})\cdots (L_r-u_m).\cr 
\endaligned
$$
Then, noting from (1.10) that
$w_{\bka}=w_{\bka_k}=w_{\bka_\dashv} s_{r, r-a_{k-1}}$, we have
 $$
\aligned
v_\bka=&v_{\bka_\dashv}(L_r-u_{k+1})\cdots (L_r-u_m)\cr 
&T_{r, r-a_1}(L_{r-a_1}-u_1)\cdots T_{r-a_{k-2}, r-a_{k-1}} 
(L_{r-a_{k-1}}-u_{k-1}).\cr
\endaligned$$
Thus, applying (3.2) to $[0,r-a_{k-1},\cdots,r-a_1,r,\cdots,r]$,  we have  
$$
v_{\bka}=v_{\bka_\dashv}(q^{c_2} L_r^{m-1} T_{r-a_{k-1},  r}^{-1}+h_1)
$$ 
where $c_2\in \Bbb Z$, $h_1\in U_1$ (see the proof of (3.2)).
On the other hand, by induction, we have, for some integer $c_1$,
$v_{\bka_\dashv}=q^{c_1} \prod_{j=1}^{r-1}L_j^{m-1} h_{(\bka_\dashv)'}^{-1}+h,$
where $h$ is a linear combinations of elements $L_1^{d_1}\cdots L_{r-1}^{d_{r-1}}
T_w$ with $w\in \frak S_{r-1}$.
Since $L_r$ commutes with any elements in $\bH_m^{r-1}$, we have 
$pr_h(q^{c_1} \prod_{j=1}^{r-1}L_j^{m-1} h_{(\bka_\dashv)'}^{-1} h_1)=0$,
$ pr_h( h q^{c_2} L_r^{m-1} T_{r, r-a_k}^{-1})=0$
and $pr_h(hh_1)=0$. Therefore, (3.5) holds (cf. (1.11)(a)).

Now we are ready to prove that the set $\{v_\bka T_w \mid w\in \fS_r\}$
is linearly independent. Suppose $\sum_{w\in \frak S_r} c_w v_{\boldkey a} 
T_w=0$.
Then, by (3.5), there is a $c\in \Bbb Z$ depending on $\bka$ such that 
$$pr_h(\sum_{w\in \frak S_r} c_w v_{\boldkey a} T_w)=
\sum_{w\in \frak S_r} c_w q^c L_1^{m-1}\cdots L_r^{m-1}
 h_{\boldkey a'}^{-1}  T_w=0.
$$ 
Because $h_{\boldkey a'}^{-1}$ and $q$ are invertible, we have  
$c_w =0 $ for every $w\in \frak S_r$ by (2.2). 
Therefore, $\{v_{\boldkey a} T_w\mid w\in \frak S_r\}$ is a linearly
independent set.  By (3.1)(d), $v_\bka\bH=v_\bka\sH$. Thus, 
 $v_{\boldkey a} \bH$ is a free $R$-module
with basis $\{v_{\boldkey a} T_w\mid w\in \frak S_r\}$.
\qed\enddemo

\proclaim{(3.6) Proposition} Let $\boldkey a\in \La[m, r] $ and write
 $v_{\boldkey a} h_{ \boldkey a'} v_{\boldkey a} =v_{\boldkey a} z_{\bka'}
=z_{\bka}v_\bka$, where $z_\bka\in\sH ({\fS_\bka})$ and $z_{\bka'}
\in\sH ({\fS_{\bka'}})$. Then
$z_\bka$ (resp. $z_{\bka'}$) is in the center of the algebra 
$\sH(\frak S_{\boldkey  a})$ (resp. $\sH(\frak S_{\boldkey  a'})$) \endproclaim

\demo{Proof} The existence of $z_\bka$ and $z_{\bka'}$ follows from (3.1a,c). 
We need only to prove that $z_{\bka'}$
is in the center of $\sH(\frak S_{\boldkey a'})$.
One may prove that $z_\bka$ is in the center of
 $\sH(\frak S_{ \boldkey a})$, similarly.

Let $T_i\in \sH(\frak S_{\boldkey a'})$. 
Then $i\neq a_m-a_j$ for  $1\le j\le m$. By (3.1)(a),  
$v_{\boldkey a} h_{ \boldkey a'} v_{\boldkey a} T_i=v_{\boldkey a} h_{ \boldkey a'}
 T_{(i) w_{\boldkey a}^{-1}} v_{\boldkey a} 
=v_{\boldkey a} T_i h_{\boldkey a'} v_{\boldkey a}=T_{(i)w_{\boldkey a}^{-1}}
 v_{\boldkey a} z_{\bka'} = v_{\boldkey a} T_i z_{\bka'}$. 
Therefore, 
$v_{\boldkey a} z_{\bka'} T_i=v_{\boldkey a}T_i z_{\bka'}$. By (3.4), 
$z_{\bka'} T_i=T_i z_{\bka'}$. So,  $z_{\bka'}$ is in the center of $\sH(\frak S_{ \boldkey a'})$.\qed\enddemo

\proclaim{(3.7) Corollary} We have 
  $T_{w_{\boldkey a}} z_{\boldkey a'} =z_{\boldkey a} 
T_{w_{\boldkey a}}$. 
\endproclaim
\demo{Proof} By (1.8), we may
write $T_{w_{\boldkey a}} z_{\boldkey a'}=zT_{w_{\boldkey a}}$
for some $z\in\sH$.
Then $v_\bka z_{\bka'}=zv_\bka$, and hence $z=z_\bka$ by (3.4).\qed
\enddemo

 \proclaim{(3.8) Corollary} For  $\boldkey a\in \La[m, r]$, 
 $v_{\boldkey a} \bH v_{\boldkey a} =v_{\boldkey a} z_{\boldkey a'} 
\sH(\frak S_{\boldkey a'})=\sH(\frak S_{\boldkey a})z_\bka v_\bka$.
\endproclaim

\demo{Proof}  Let $\iota$ be the  anti-automorphism of $\bH$ defined in (2.4).
By (3.1a), we have $\tpi_{\boldkey a'} \sH\pi_{\boldkey a} =
\iota (\pi_{\boldkey a}\sH \tpi_{\boldkey a'})=\iota(\sH(\frak S_{\bka}) v_\bka)
=\tpi_{\boldkey a'} h_{\boldkey a'} \pi_{\boldkey a}\sH(\frak S_{\boldkey a})$.
So, 
$$
\aligned v_{\boldkey a} \bH v_{\boldkey a}& =v_{\boldkey a}\sH v_{\boldkey a}
=\pi_{\boldkey a} h_{\boldkey a} \left( \tpi_{\boldkey a'} \sH 
\pi_{\boldkey a}\right) h_{\boldkey a} \tpi_{\boldkey a'}\cr
&=\pi_{\boldkey a} h_{\boldkey a} \tpi_{\boldkey a'} 
h_{\boldkey a'} \pi_{\boldkey a} \sH(\frak S_{\boldkey a}) h_{\boldkey a} \tpi_{\boldkey a'}\cr
&= v_{\boldkey a} h_{\boldkey a'} v_{\boldkey a} \sH(\frak S_{\boldkey a'})
=v_{\boldkey a} z_{\boldkey a'} \sH(\frak S_{\boldkey a'}).\qed\cr
\endaligned
$$
\enddemo

\proclaim{(3.9) Proposition}
 Let $\boldkey a\in \La[m, r]$. Then the following are equivalent:

(a) $z_{\boldkey a}$
(equivalently $z_{\boldkey a'}$) is invertible.

(b) $e_{\boldkey a} =v_{\boldkey a} h_{\boldkey a'} z_{\boldkey a}^{-1}$
is an idempotent.

(c) $v_{\boldkey a} \bH$ is a projective right
 $\bH$-module.\endproclaim

\demo{Proof}
Let $J=v_{\boldkey a} \bH$. By (3.1)(d), we have $J^2=v_{\boldkey a} z_{\boldkey a'} \sH$. 
By (3.4), the map
from $\sH$ to $v_{\boldkey a} \sH$ sending  
$h\rightarrow v_{\boldkey a} h$ for  $h\in \sH$ is  injective. 
Therefore, $J^2=J$ if and only if 
$z_{\boldkey a'} \sH=\sH$, which is equivalent to say
 that $z_{\boldkey a'}$ is invertible. 

If $J$ is projective, then $J=e\bH$ for some idempotent $e\in \bH$, 
and $J=J^2$. So, $z_{\boldkey a'}$,  and hence $z_{\bka}$, is invertible. 
Now, assume $z_{\boldkey a}$ is 
invertible. Putting $e_{\boldkey a} =v_{\boldkey a} h_{\boldkey a'} z_{\boldkey a}^{-1}$ and recalling $h_{\boldkey a'}=T_{w_{\bka'}}$,
we have 
$$\aligned 
e_{\boldkey a}^2 & =v_{\boldkey a} h_{\boldkey a'} z_{\boldkey a}^{-1} 
v_{\boldkey a} h_{\boldkey a'} z_{\boldkey a}^{-1}=v_{\boldkey a} h_{\boldkey a'} z_{\boldkey a}^{-1} 
\pi_{\boldkey a} h_{\boldkey a} \tpi_{\boldkey a' } h_{\boldkey a'}
 z_{\boldkey a}^{-1}\cr
&=v_{\boldkey a} h_{\boldkey a'} \pi_{\boldkey a} z_{\boldkey a}^{-1} h_{\boldkey a} \tpi_{\boldkey a'} 
h_{\boldkey a'} z_{\boldkey a}^{-1}\cr
&=v_{\boldkey a} h_{\boldkey a'} v_{\boldkey a} z_{\boldkey a'}^{-1} h_{\boldkey a'} z_{\boldkey a}^{-1}
=v_{\boldkey a} h_{\boldkey a'} z_{\boldkey a}^{-1}=e_{\boldkey a}.\cr\endaligned
$$  
 Thus, 
$e_{\boldkey a}$ is an idempotent and $J=v_{\boldkey a} 
\bH=e_{\boldkey a} \bH$. Therefore, $v_{\boldkey a}\bH$ is projective.
\qed\enddemo

\proclaim{(3.10) Proposition} Let  $\boldkey a, \boldkey b\in \La[m, r]$
and assume that $z_\bka, z_\bkb$ are invertible.
 Then $e_{\boldkey a}\bH e_{\boldkey b}=0$
unless $\boldkey a=\boldkey b$. 
Moreover, we have 
$e_{\boldkey a}\bH e_{\boldkey a} \cong \sH(\frak S_{\boldkey a})$.\endproclaim

\demo{Proof} 
By (2.8) and (3.1)(c),  $v_{\boldkey a}\bH v_{\boldkey b}\neq 0$ implies
$\boldkey a\lecq \boldkey b$ and  $\boldkey a'\lecq \boldkey b'$. So 
   $\boldkey a= \boldkey b$.
The second assertion follows from   the following equality:
$$
\split e_{\boldkey a} \bH e_{\boldkey a} &  
           =v_{\boldkey a} \bH v_{\boldkey a} h_{\boldkey  a'} {
 z_{\boldkey a}}^{-1}
           =  v_{\boldkey a}   \sH(\frak S_{ \boldkey a'} ) h_{\boldkey a'}\text{ by (3.8) }
             \cr
           &=v_{\boldkey a} h_{\boldkey a'}\sH(\frak S_{\boldkey a} )
            \cong \sH(\frak S_{\boldkey a}).\qed \cr\endsplit
$$ 
\enddemo

\proclaim{(3.11) Corollary} Assume that
$z_\bka$ are invertible for all $\bka\in\La[m,r]$ and
let $\varepsilon=\sum_{\boldkey a\in \La[m, r]} 
e_{\boldkey a}$. Then $\varepsilon \bH\varepsilon\cong
\oplus_{\bka\in \La[m, r]} \sH(\frak S_{\bka})$. 
\endproclaim 
\demo{Proof}  The result  follows from (3.10).
\qed\enddemo

\subhead 4. The Poincar\'e polynomial of $W$\endsubhead 
To generalize the Morita equivalence (1) given
in the introduction, we need two more ingredients:
First,  we want to know when the hypothesis in (3.11) holds.
This leads to the introduction of the Poincar\'e polynomial of the bottom.
Second,  we need to prove that the direct sum $\oplus_{\bka\in \La[m, r]} 
v_\bka \bH$ is a projective generator for the category $\bH\hmod$.
The latter requires that $R$ is an integral domain.

\proclaim{(4.1) Lemma}  For  $\boldkey r_i=(0,
\underset {i-1} \to { \underbrace {0\cdots 0}}, r \cdots , r)\in \La[m, r]$, 
let  $z_{\boldkey r_i'}\in \frak S_{\boldkey r_i'}$ defined in (3.6). Then 
$z_{\boldkey r_i'}= \prod\Sb 1\le j\le m \\ j\neq  i\endSb 
 \prod_{k=1}^r  (u_i q^{1-k} T_{k, 1} T_{1, k}-u_j)$.
\endproclaim

\demo{Proof} Since $w_{\boldkey r_i}=1$,
we have $h_{\boldkey r_i}=T_{w_{\boldkey r_i}}=1$.
 It follows from (2.1)(5) and (2.7)(a) that 
$$v_{\boldkey r_i}^2 =v_{\boldkey r_i} 
\prod\Sb 1\le j\le m\\ j\neq i\endSb 
\prod_{k=1}^r (u_i q^{1-k} T_{k, 1} T_{1, k}-u_j).
$$
On the other hand, by (3.6), we have 
$v_{\boldkey r_i}^2 =v_{\boldkey r_i} z_{\boldkey r_i'}$.
Now, the result follows from (3.4) immediately.
\qed\enddemo

\definition{(4.2) Definition} For positive integers $m, r$ and
$i=1,\cdots,m$, let 
 $$\aligned 
&f_{m, r, i}= f_{m, r, i} (u_1, \cdots u_m, q)= \prod\Sb  1\le j\le m\\ i\neq j\endSb  
  \prod_{k=1-r}^{r-1} (u_i q^k-u_j),\cr
  &  f_{m, r}= f_{m, r} (u_1, \cdots, u_m, q)
=\prod_{i=1}^{m-1} \prod_{j=i+1}^m  \prod_{k=1-r}^{r-1}
   (u_i q^k-u_j).\cr\endaligned
  $$ 
We shall see below that the polynomial $f_{m,r}$ 
can be viewed as the Poincar\'e polynomial of the bottom $C$.
Let $d_{\fS_r}$ be the Poincar\'e polynomial of $\fS_r$, i.e.,
$d_{\fS_r}=\sum_{w\in\fS_r}q^{l(w)}$, then the polynomial
$d_W=f_{m,r}d_{\fS_r}$ is called the {\it Poincar\'e polynomial} of
the complex reflection group $W$.
\enddefinition

\proclaim{(4.3) Proposition} Maintain the notation introduced above.  
The element $z_{\boldkey r_i'}$ is invertible if and only if 
$f_{m, r, i}$ is invertible in $R$.\endproclaim

\demo{Proof} 
By (4.1), we see that the invertiblity of $z_{\boldkey r_i'}$ is 
equivalent to the invertiblity of  
$h_{ij}=\prod_{k=1}^r (u_iq^{1-k} T_{k,1}T_{1,k}-u_j)$ for 
all $j\neq i$.  By \cite{DJ2, (4.3)}, we see that,
if $R$ is a field, then $h_{ij}$ in invertible if and only if 
$f_{m, r,i}$ is invertible. The general case
follows by an argument similar to the one for \cite{DJ2, (4.5)}.   
\qed\enddemo

We need some preparation in order to get the main result of this paper.
 As before, we assume that $R$ is a commutative ring (with 1).
For any $\bka\in\La[m,r]$ and $x_1,\cdots,x_{m-1}\in R$,
recall the element defined after (2.6)
$$\pi(\bka)=\pi(\bka;x_1,\cdots,x_{m-1})=\pi_{a_1}(x_1)\cdots
\pi_{a_{m-1}}(x_{m-1}).$$
Write $\pi(\bka)$ as a polynomial in $L_i$. The degree of this polynomial
is denoted $\deg_i(\pi(\bka))$.

\proclaim{(4.4) Lemma} Maintain the notation above and assume 
$x_1,\cdots,x_{m}$ are  a permutation of $u_1,\cdots,u_m$.
For  $\boldkey a \in \La [m, r]$, the right 
ideal $\pi({\boldkey a}) \bH$ is spanned by      
$\sB_{\boldkey a} = \{\pi({\boldkey a}) L_1^{c_1}\cdots L_r^{c_r} T_w\mid 
 w\in \frak S_r,\deg_j \pi({\boldkey a}) 
 +c_j\le m-1, \forall j\}$. 
\endproclaim

\demo{Proof}
Let $M_\bka$ be the $R$-submodule spanned by $\sB_\bka$. 
Then $M_\bka\subseteq \pi({\boldkey a})\bH$. Since $\pi(\bka)\in M_\bka$,
it suffices to prove that 
$M_\bka$ is a right $\bH$-module, or equivalently
$M_\bka T_0\subseteq M_\bka$.  

For  $ \pi({\boldkey a}) L_1^{c_1}\cdots L_r^{c_r} T_w\in \sB_\bka$. We write
$w=s_{i, 1} y$ with $y\in \frak S_{\{2, \cdots, r\}}$ (cf. (1.3)). 
Since $M_\bka$
is a right $\sH$-module and $T_0T_y=T_yT_0$, we  only need to  prove 
$\pi({\boldkey a}) L_1^{c_1}\cdots L_r^{c_r} T_{i, 1}T_0\in M_{\bka}$, which
is equivalent to  
$$\pi({\boldkey a}) L_1^{c_1}\cdots L_r^{c_r} L_i\in M_\bka \text{ 
for every $1\le i\le r$}. \tag 4.5
$$
  If $a_{m-1}=0$, then $\pi(\bka)=1$ and 
$\pi(\bka)\bH=\bH$.  By (2.2), 
(4.5) holds in this case.  We now assume 
$a_{m-1}\neq 0$.
Apply induction on $i$. The last relation in (2.1) implies (4.5) for $i=1$.
(If $\deg_1(\pi(\bka))+c_1+1=m$, then $\pi(\bka)$ contains part of
the product in (2.1)(5) involving parameters $x_i,\cdots,x_{m-1}$. 
Write $L_1^{c_1+1}=\prod_{k}(L_1-u_{i_k})+\sum_{j\le
c_1}\al_j L_1^j$,
where $\{u_{i_k}\}=\{x_1,\cdots, x_{i-1},x_m\}$. Now, (2.1)(5) implies (4.5)
in this case.)

We assume now that $i>1$ and (4.5) holds for all $L_j$ with $j<i$.
The case for
 $\deg_i \pi({\boldkey a})  +c_i<m-1$ is trivial. 
Suppose $\deg_i \pi({\boldkey a})  +c_i=m-1$.
Let $k$ be  the  integer 
with  $a_{k-1}<i\le a_{k}$. 
Then, $\deg_i \pi({\boldkey a}) =m-k$ and  
$\pi_{a_j}(x_j) T_{i, 1}=T_{i, 1} \pi_{a_j}(x_j) $ 
for all $j\ge k$. On the other hand, by (2.5)(d) and induction on $i$, 
we have for some $c$
$$
L_i^l- q^c T_{i, 1} L_1^l T_{1, i}\in U_{i, l}\tag 4.6
$$
where  $U_{i, l} $ is the   
free $R$-submodule  spanned by 
$\{L_1^{d_1}\cdots L_i^{d_i}T_w\mid 0\le d_j<l , \forall j\le i, 
 w\in \frak S_i\}$. 
In particular, 
$L_{i}^{c_i+1}=q^cT_{i, 1} L_1^{c_i+1}T_{1, i}+h$ for some $c$ and 
$h\in U_{i, c_{i}+1}$. So we have   
$$
\aligned 
&\pi(\bka) L_1^{c_1} \cdots L_{i-1}^{c_{i-1}} L_i^{c_i} L_{i+1}^{c_{i+1}}\cdots L_r^{c_r}L_i\cr 
=&\pi(\bka) L_1^{c_1} \cdots L_{i-1}^{c_{i-1}} L_i^{c_i+1} L_{i+1}^{c_{i+1}}\cdots L_r^{c_r}\cr
=&\pi(\bka) L_1^{c_1} \cdots L_{i-1}^{c_{i-1}}(q^cT_{i, 1} L_1^{c_i+1}T_{1, i}+h) L_{i+1}^{c_{i+1}}\cdots L_r^{c_r}\cr
=&q^c\pi(\bka) L_1^{c_1} \cdots L_{i-1}^{c_{i-1}}T_{i, 1} L_1^{c_i+1}T_{1, i} L_{i+1}^{c_{i+1}}\cdots L_r^{c_r}\cr
&+\pi(\bka) L_1^{c_1} \cdots L_{i-1}^{c_{i-1}}h L_{i+1}^{c_{i+1}}\cdots L_r^{c_r}\cr
=&X_1+X_2\cr
\endaligned 
$$
%
%= &  q^{1-i} \pi(\bka) \left( \prod_{j=1}^{i-1} L_{j}^{c_j}\right) T_{i, 1} L_1^{c_i+1} T_{i, 1}\prod_{j=i+1}^r 
% L_{j}^{c_{j}}+ \pi(\bka)\left( \prod_{j=1}^{i-1}  L_j^{c_j}\right)   h \prod_{j=i+1}^r   L_{j}^{c_{j}}\cr
%=& q^{1-i}\prod_{j=1}^{i-1} L_j^{c_j}  \prod_{j=1}^{k-1} \pi_{a_j}(x_j)
%\left( T_{i, 1}
%\prod_{j=k}^{m-1} \pi_{a_j}(x_j) L_1^{c_i+1} T_{i, 1}\right)
% \prod_{j=i+1}^r  L_{j}^{c_j }\cr
%+ &\pi(\bka) \prod_{j=1}^{i-1} L_j^{c_j}  h \prod_{j=i+1}^r  L_j^{c_j}\cr 
%
We now prove that each of the last two terms $X_i$
 above is in $M_{\bka}$. 
For $h=L_1^{d_1}\cdots L_i^{d_i} T_w\in U_{i, c_i+1}$, we have 
$T_w L_j=L_jT_w$ for all $j\ge i+1$ and $d_i<c_i+1$, and so, 
 $\deg_i \pi(\bka)+ d_i <\deg_i \pi(\bka)+c_i+1=m$. 
On the other hand, 
we have, by definition, $\deg_j \pi(\bka)+c_j\le m-1$ for all $1\le j\le m$.
Therefore, by inductive hypothesis,
$$
\pi(\bka) L_1^{c_1} \cdots L_{i-1}^{c_{i-1}}(L_1^{d_1}\cdots L_i^{d_i} T_w) 
 L_{i+1}^{c_{i+1}}\cdots L_r^{c_r}
=\pi({\bka})\prod_{j<i} L_j^{c_j+d_j} L_i^{d_i} \prod_{j>i} L_j^{c_j}  T_w\in M_\bka.
$$
This proves $X_2\in M_\bka$.

Since $i\le a_k$, we have

$$\aligned
X_1&=
\prod_{j=1}^{i-1} L_j^{c_j}  \prod_{j=1}^{k-1} \pi_{a_j}(x_j)
\left( T_{i, 1}
\prod_{j=k}^{m-1} \pi_{a_j}(x_j) L_1^{c_i+1} T_{i, 1}\right)
 \prod_{j=i+1}^r  L_{j}^{c_j }\cr
&=X_{11}(X_{12})X_{13}.\cr
\endaligned
$$
where $X_{11}$ (resp. $X_{13}$) denotes the first two  products 
(resp. last  product) and $X_{12}$ denotes the element in the parenthesis.
As in the argument for $i=1$,   
$L_1^{c_i+1 }\prod_{j=k}^{m-1} \pi_{a_j}(x_j)$ is
 an $R$-linear combinations of $ L_1^{l} \prod_{j=k}^{m-1} \pi_{a_j}(x_j)$ 
with $l<c_i+1 $, while, by (4.6),
$L_i^l-T_{i,1}L_1^lT_{1,i}=h\in U_{i,l}$.
Thus,   
$T_{i, 1}  L_1^{l}\prod_{j=k}^{m-1} \pi_{a_j}(x_j) T_{i, 1}=
\prod_{j=k}^{m-1}    \pi_{a_j}(x_j)L_i^l+  
\prod_{j=k}^{m-1}  \pi_{a_j}(x_j)h$. Hence
$X_1$  can be expressed as   a linear combination of 
$X_{11}( \prod_{j=k}^{m-1}
  \pi_{a_j}(x_j)L_i^l)X_{13}$
and $X_{11}( \prod_{j=k}^{m-1}\pi_{a_j}(x_j)h)X_{13}$, where $l\le c_i$ and
$h\in U_{i,l}$.
Since $l\le c_i$, we have clearly $X_{11}( \prod_{j=k}^{m-1}
  \pi_{a_j}(x_j)L_i^l)X_{13}\in
M_\bka$. By inductive hypothesis, we have
$X_{11}(\pi_{a_j}(x_j)h)X_{13}\in M_\bka$, too.
Therefore, $X_1\in M_\bka$, proving (4.5).
\qed\enddemo

For   $\boldkey a=[a_i] 
\in \La[m, r]$, 
let ${\boldkey a_i}$
be  defined as in (1.9). Recall $h_\bka=T_{w_\bka}$.

\proclaim{(4.7) Proposition} For $\bka\in\La[m,r]$, 
write $\bka_\dashv=[b_0,\cdots,b_m]$ as in (1.10).   
Let  $V_i= \pi_{\boldkey a_\dashv} h_{\boldkey a_\dashv} T_{r, r-b_{i-1}} \tpi_{ \boldkey a_i'}\bH$.   

Then

(a) $V_{1}=v_{\boldkey a_\dashv}\bH$ and $V_m=v_{\bka_m}\bH
$. 

(b) The set 
$\sB_i=\{\pi_{\boldkey a_\dashv} h_{\boldkey a_\dashv} T_{r, r-b_{i-1}} \tpi_{\boldkey a_i'}
L_{r-b_{i-1}}^c T_w\mid c\le m-i, w\in \frak S_r\}$ is a basis of $V_i$. 
 In particular,
the rank  of $V_i$ is $(m-i+1)r! $.

(c) $V_{i+1}$ is a pure $R$-submodule of $V_i$.
\endproclaim

\demo{Proof} We first note that
$b_i=a_i, 1\le i\le k-1, b_i=r-1, k\le i\le m,$
where $k$ is the minimal index with $a_k=r$. Note also from (1.9) and (1.2)
that 
  $$\aligned 
\bka_1'&=[0, \cdots, 0, r-a_{k-1}-1, \cdots, r-a_{1}-1,  r] \cr
(\bka_\dashv)'&=[0, \cdots, 0, r-a_{k-1}-1, \cdots, r-a_{1}-1, r-1].\cr
\endaligned
$$
Thus,  $\tpi_{\bka_1'}=\tpi_{(\bka_\dashv)'}$ by definition (2.6), and
$V_1=v_{\bka_\dashv}\bH$.
Because $h_{\boldkey a_\dashv} T_{r, r-b_{m-1}}=h_{\bka_m}$ (see (1.10))
and $\pi_{\bka_m}=\pi_{\bka_\dashv}$,   
 we have  $V_m=v_{\bka_m}\bH$, proving (a).

For (b), we first prove that $\sB_i$ spans $V_i$.
Let $V_i'$ be the submodule spanned by $\sB_i$.
Applying (4.4) to $\tpi_{ \boldkey a_i'}\bH$, it suffices to prove  
$$\pi_{\boldkey a_\dashv} h_{\boldkey a_\dashv} T_{r,  r-b_{i-1} } 
\tpi_{ \boldkey a_i'}L_1^{c_1}\cdots L_{r-b_{i-1}}^{c_{r-b_{i-1}}}
\cdots L_r^{c_r}\in V_i',
\tag4.8$$
where $c_j\le m-1-\deg_j\tpi_{\bka_i'}$, $\forall j$. Since
$\deg_{r-b_{i-1}}\tpi_{\bka_i'}=i-1$, we have, in particular, 
$c_{r-b_{i-1}}\le m-i$.

We first look at the elements 
$\pi_{\boldkey a_\dashv} h_{\boldkey a_\dashv} T_{r,  r-b_{i-1} } 
\tpi_{ \boldkey a_i'}L_j$ with $j\neq r-b_{i-1}$.
Suppose $j< r-b_{i-1}$ and write
$T_{r, r-b_{i-1}}\tpi_{\bka_i'}=
T_{r, r-b_{i-1}}\tpi_{(\bka_{\dashv})'}(L_{r-b_1}-u_1) \cdots 
(L_{r-b_{i-1}}-u_{i-1})=
\tpi_{(\bka_\dashv)'} h_i$, where $ h_i=T_{r, r-b_1} (L_{r-b_1}-u_1) \cdots 
T_{r-b_{i-2}, r-b_{i-1}} (L_{r-b_{i-1}}-u_{i-1})$. 
Then $L_jh_i=h_i L_j$ and $h_i T_l=T_l h_i$ for all
$s_l\in \frak S_j$. Thus,
  by (3.1)(b),    
$$
\pi_{\boldkey a_\dashv} h_{\boldkey a_\dashv} T_{r,  r-b_{i-1} } 
\tpi_{ \boldkey a_i'}L_j
=v_{\bka_\dashv} L_j h_i\in 
\pi_{\boldkey a_\dashv} h_{\boldkey a_\dashv} T_{r,  r-b_{i-1} } \tpi_{ \boldkey a_i'} 
\sH(\frak S_j).
\tag 4.9
$$
Suppose now $j\ge r-b_{i-1}+1$. Choose $l$ such that
$r-b_{l}+1 \le j\le r-b_{l-1}$. Then $l\le i-1$ and $l$ is the minimal index
having property $b_l=b_{i-1}$.
Using definition (2.6), we have $\tpi_{\bka'_i}(L_{r-b_l+1}-u_l)=
\tpi_{\boldkey c'}$, where
$\boldkey c=[0, b_1, \cdots, b_{l-1}, b_l-1, b_{l+1},
\cdots,b_m]+{\boldkey 1}_i \in \La[m, r]$ (see (1.9)).
Since $\pi_{\bka_\dashv}=\pi_{\bka_m}$ and $\bka_m\not\lecq\bkc$ (as $l<i$),
 we have, by (2.8) that 
$\pi_{\boldkey a_\dashv} h_{\boldkey a_\dashv} T_{r,  r-b_{i-1} } 
\tpi_{\boldkey a_i'}(L_{r-b_l+1}-u_{l})=0$. 
Therefore, for some $c$,  
$$
\pi_{\boldkey a_\dashv} h_{\boldkey a_\dashv} T_{r,  r-b_{i-1} } \tpi_{\boldkey a_i'} L_j  
=u_{l}q^{c} \pi_{\boldkey a_\dashv} h_{\boldkey a_\dashv} T_{r,  r-b_{i-1}} 
\tpi_{ \boldkey a_i'} 
T_{j, r-b_{l}+1} T_{r-b_{l}+1, j}.  
\tag 4.10
$$
Now, noting (2.5)(a), we see that (4.8) follows from (4.9) and (4.10),
since the elements in $\sH(\fS_j)$ occurring in (4.8) and the element
$T_{j, r-b_{l}+1} T_{r-b_{l}+1, j}$ occurring in (4.9) commute with $L_{j+1}$
and $L_{r-b_{i-1}}$.

To see the linear independence, we apply induction on $i$. For $i=1$, suppose 
$$
\sum\Sb 0\le c\le m-1\\ w\in \frak S_r \endSb 
 f_{c, w} v_{\boldkey a_\dashv} L_r^c T_w=0.
$$ 
 Since  $pr_h (v_{\boldkey a_\dashv})= q^d 
\prod_{j=1}^{r-1} L_j^{m-1} h^{-1}_{\boldkey (a_\dashv)' }$ for some $d\in \Bbb Z$
depending on $\bka_\dashv$ (see (3.5)), we have  
$$
pr_{\bkc} (\sum\Sb 0\le c\le m-1\\ w\in \frak S_r \endSb 
 f_{c, w} v_{\boldkey a_\dashv} L_r^c T_w)=
q^d\sum_{ w\in \frak S_r} 
f_{c, w} L_1^{m-1}\cdots L_{r-1}^{m-1}L_r^c h_{(\bka_\dashv)'}^{-1} T_w=0,
$$ 
where $\bkc=(m-1, \cdots, m-1, c)$.
Thus, $f_{c, w}=0$ for all $c$ and $w$,
 since both $q$ and $h^{-1}_{\boldkey (\bka_\dashv)' }$ 
are invertible. Therefore,  $\sB_1$ is a basis of $V_1$.

 We assume now that 
$\sB_i$ is a linearly independent set.  We hope to prove that the set  
 $\sB_{i+1}$ is linearly independent, too.
Using induction on $l$, we  can easily  prove   
$$
T_{a, b} L_b^l =q^{a-b} L_a^l T_{b, a}^{-1} + \text{ terms involving }
L_a^{c_a}\cdots L_{b}^{c_b}T_w\text{ for $a\ge b$ }
$$ 
where $w\in \frak S_{\{b, b+1, \cdots, a\}}$ and all $c_i<l$.
Using this and noting 
 $\tpi_{\bka_{i+1}'}=\tpi_{\bka_i}' (L_{r-b_{i}}-u_i)$,  we have 
$$
\aligned &  \pi_{\boldkey a_\dashv} h_{\boldkey a_\dashv} T_{r, r-b_i} 
\tpi_{\boldkey a_{i+1}'} L_{r-b_i}^c T_w\cr  
= &  q^{a_i-a_{i-1}} \pi_{\boldkey a_\dashv} h_{\boldkey a_\dashv} T_{r, r-b_{i-1}} 
\tpi_{\boldkey a_i'} L_{r-b_{i-1}}^{c+1}
T_{r-b_i, r-b_{i-1}}^{-1}  T_w + \cr 
&+ \text{ terms involving }  \pi_{\boldkey a_\dashv} h_{\boldkey a_\dashv} T_{r, r-b_{i-1}} 
\tpi_{\boldkey a_i'}
 L_{r-b_{i-1}}^{l} T_z,
\cr\endaligned \tag 4.11
$$
where $l<c+1$ and $z\in \frak S_r$.
By induction, $\sB_i$ is a basis of $V_i$. Because $T_{r-b_i, r-b_{i-1}}$ is invertible, 
the set $\{\pi_{\boldkey a_\dashv} h_{\boldkey a_\dashv} T_{r, r-b_{i-1}} \tpi_{\boldkey a_i'}
L_{r-b_{i-1}}^c T_{r-b_i, r-b_{i-1}}^{-1} T_w\mid c\le m-i, w\in \frak S_r\}$ is linearly independent in $V_i$. Therefore, the condition
$l<c+1$ in (4.11) implies that   $\sB_{i+1}$ 
is linear independent. This completes the proof of (4.7)(b).

We now prove the purity in (c).
Since  $$\pi_{\boldkey a_\dashv} h_{\boldkey a_\dashv} T_{r, r-b_i} \tpi_{\boldkey a_{i+1}'}
=\pi_{\boldkey a_\dashv} h_{\boldkey a_\dashv} T_{r, r-b_{i-1}}\tpi_{\boldkey a_i' } 
T_{r-b_{i-1}, r-b_i}
(L_{r-b_i}-u_i),$$ 
$V_{i+1}$ is a submodule of $V_i$. Write 
$h=\sum f_{c,w} \pi_{\boldkey a_\dashv } h_{\boldkey a_\dashv} T_{r, r-b_{i}} 
\tpi_{\boldkey a_{i+1}'}
L_{r-b_{i}}^c T_w$ for $h\in V_{i+1}$. 
To show that $V_{i+1}$ is a pure $R$-submodule, we need to show 
that $x\in R$ divides any of 
$f_{c,w}$ if $h\in xV_{i}$. 
We prove it by induction on $c_0$,  the highest degree of 
$L_{r-b_i}$ in the expression of 
$h$. Since $\sB_i$ is a basis of $V_i$  and 
$q$ is invertible,
 $x$ divides $f_{c_0,w}$ by (4.11).
Thus,  $h-\sum f_{w, c_0} \pi_{\boldkey b} h_{\boldkey b} T_{r, r-b_i}
 \tpi_{\boldkey a_{i+1}'} L_{r-b_i}^c T_w \in xV_{i}$.
Now, the result  follows from induction.
\qed\enddemo

From here onwards, we assume that $R$ is an {\it integral domain}.

\proclaim{(4.12) Lemma} Keep the notation above. 
 We have short  exact sequence
$$
0\rightarrow V_{i+1}\rightarrow V_i\rightarrow v_{\boldkey a_i}\bH\rightarrow 
0.
$$
\endproclaim

\demo{Proof} Let
$$
h_i=(L_{b_{i}+1}-u_{i+1}) 
T_{b_{i}+1, b_{i+1}+1}\cdots (L_{b_{k-1}+1}-u_k) T_{b_{k-1}+1, r}
\prod_{l=k+1}^m   (L_r-u_{l}),
$$
where, as usual, $k$ is the minimal index such that $a_k=r$. Then,
 by (2.5)(a)-(c), we have
$ h_i \pi_{\bka_\dashv}=\pi_{\bka_i} T_{b_i+1, r} $.
Define  $\psi_i: V_i\rightarrow v_{\boldkey a_i}\bH$
by setting $\psi_i(h)=h_i h$, $h\in V_i$.
Note that $\psi_i$ is well-defined by (1.10), and clearly,
$\psi_i$ is  surjective.
 Because  $\boldkey a_{i}\rec   \boldkey a_{i+1}$, we have  
$V_{i+1}\subseteq \ker \psi_i$ by (2.8) and (3.1)(c).
By (3.4) and (4.7)(b), we have 
$\dim V_i=\dim V_{i+1} +\dim v_{\boldkey a_i}\bH$, 
forcing  $V_{i+1}=\ker{\psi_i}$ over the quotient field $F$ of  $R$.
 However, by (4.7)(c), $V_{i+1}$ is a pure $R$-submodule of $V_i$. 
Thus $V_{i+1}=\ker{\psi_i}$, proving  the required short exact sequence.
\qed
\enddemo

\proclaim{(4.13) Lemma} Let $\bka=[a_i]\in \La[m, r]$ and 
write $\bka_\dashv=[b_i]$. Assume 
$k$ is the minimal index with $a_k=r$. Define     
$$\tilde V_{i}=\pi_{\bka_\dashv} h_{\bka_\dashv} T_{r, r-b_{i-1} } \tpi_{\bka_{i}'}
\prod_{j=i+1}^m \pi_{r-b_{i-1}}(u_j) \bH, \text{ for }
i=k,\cdots, m-1.
$$
If  $f_{m, r}$ is a unit in the 
integral domain $R$, then  we have $V_{i}= V_{i+1}\oplus \tilde V_i$
and $\tilde V_i\cong v_{\bka_i}\bH$ for 
all $i=k,\cdots, m-1$.
\endproclaim

\demo{Proof} Let $\psi_{i} : V_i\rightarrow v_{\bka_i}\bH, i\ge k$, be the 
$\bH$-module homomorphism defined  in the proof of (4.12).
Then, $
\psi_{i}(\tilde V_{i} )=v_{\bka_{i}}\prod_{j=i+1}^m \pi_{r-b_{i-1} } (u_j)\bH
$.
By (3.1)(b), we have $v_{\bka_i}L_1=u_iv_{\bka_i}$ for all $i\ge k$.
Also, for $j<r-a_{k-1}$, $s_j\in\fS_{\bka_k'}$. So $v_{\bka_k}T_j=
T_{j'}v_{\bka_k}$ (see (3.1)(a)), and hence, 
$v_{\bka_k}L_l=u_kq^{1-l}v_{\bka_k}T_{l,1}
T_{1,l}$ for all $l=1,\cdots,r-b_{k-1}$, noting $a_{k-1}=b_{k-1}$.
Therefore, we have
$$
\psi_{i}(\tilde V_{i} )=v_{\bka_{i}}\prod_{j=i+1}^m \pi_{r-b_{i-1} } (u_j)\bH
=v_{\bka_{i}}\prod_{j=i+1}^m  \prod_{l=1}^{r-b_{i-1}} 
(q^{1-l} u_{i} T_{l, 1}T_{1, l}-u_j)\bH.
$$
Since $\prod_{j=i+1}^m  \prod_{l=1}^{r-b_{i-1} } 
(q^{1-l} u_{i} T_{l, 1}T_{1, l}-u_j)$ is a factor of $z_{\boldkey r_i'}$ (see (4.1)),
and  $f_{m, r}\in R $ is a unit, $\prod_{l=1}^{r-b_{i-1} } 
(q^{1-l} u_{i} T_{l, 1}T_{1, l}-u_j)$ is invertible by (4.3). Therefore,  
$\psi_{i}(\tilde V_{i}) =v_{\bka_{i}}\bH$. Thus,
by (4.12), the rows of the following commutative diagram  
$$\matrix
0&\rightarrow& V_{i+1}&\rightarrow& \tilde V_i+V_{i+1}&
\overset {\psi_i}\to {\rightarrow}& v_{\bka_i}\bH&\rightarrow& 0\cr
&&\Vert&&\downarrow&&\Vert&&\cr
0&\rightarrow& V_{i+1}&\rightarrow& V_i&
\overset {\psi_i}\to {\rightarrow}& v_{\bka_i}\bH&\rightarrow& 0.
\endmatrix
$$
are exact. So the short-five lemma implies $\tilde V_i+V_{i+1}= V_i$. 
On the other hand, consider 
the element $x= \tpi_{\bka_{i}'}
\prod_{j=i+1}^m \pi_{r-b_{i-1}}(u_j)$ and
the multiplication
$ xL_{r-b_{i-1}}$.
Clearly, $x$ has a factor $(L_1-u_1)\cdots(L_1-u_{i-1})(L_1-u_{i+1})\cdots
(L_1-u_{m})$. So,
if $i>k$, then $r-b_{i-1}=1$, and hence, $xL_1=u_ix$,
while, for $i=k$, we have $xT_{r-b_{k-1},1}=T_{r-b_{k-1},1}x$ and
consequently, $xL_{r-b_{k-1}}=q^c u_i xT_{r-b_{k-1},1}T_{1,r-b_{k-1}}$.
Therefore,
by (4.7)(b), $\tilde V_{i}$ is spanned by 
$\{\pi_{\bka_\dashv} h_{\bka_\dashv} T_{r, r-b_{i-1}} \tpi_{\bka_{i}'}
\prod_{j=i+1}^m \pi_{r-b_{i-1}}(u_j)T_w\mid w\in \frak S_r\}$,
which is linearly independent, since its
 image under  the homomorphism $\psi_i$
is linearly independent. Hence,   
$\tilde V_i$ is $R$-free of rank $r!$, forcing the sum
$\tilde V_i+ V_{i+1}$ is a direct sum, i.e.,
$V_i=\tilde V_i\oplus V_{i+1}$, and consequently, 
$\tilde V_i \cong v_{\bka_i}\bH$. \qed\enddemo

We now prove  the main result of the paper. The case when $\bH$ is
the Hecke algebra of type $B$ was first obtained in \cite{DJ2}.

\proclaim{(4.14) Theorem} Let $\bH=\bH_m^r$ be the Ariki-Koike  algebra
over an integral domain $R$.  
Assume  $f_{m, r}\in R$ is a unit.

(a) For every $\boldkey a\in \La[m, r_1]$ with $r_1\le r$, the right ideal 
$v_{\boldkey a}\bH$ is projective.

(b) For any given $r_1\le r$, each projective indecomposable 
$\bH$-module is isomorphic to a 
direct summand of $v_{\boldkey a}\bH$ for some $\boldkey a\in \La[m, r_1]$.

(c) The categories of $\bH$-modules and 
$\oplus_{\la\in \La(m, r)} \sH(\frak S_{\la})$-modules
are Morita equivalent.\endproclaim

\demo{Proof} For $r_1\le r$, $\bH$ is free over $\bH^{r_1}_m$. So the
functor $-\otimes_{\bH_m^{r_1} } \bH$ is exact. Therefore,
it suffices to prove (a) for $r_1=r$.

We apply induction on $r$. For $r=1$, we have
$\bka=\boldkey r_i$ for some $i$. By the invertibility
of $f_{m,r}$, (4.3) and (3.8), we see that
$v_\bka\bH$ is projective. 
Assume now that $r>1$ and the result holds for $r-1$. We now prove the 
projectivity of $v_\bka\bH$ by downward induction on the partial ordering
$\recq$. If $\bka=\boldkey r_1=[0,r,\cdots,r]$, 
the maximal element of $\La[m,r]$,
then $v_\bka\bH$ is projective by (4.3) and (3.8).
Suppose now that $\boldkey r_1\succ\bka$ and that 
 $ v_{\bkb}\bH$ is projective for all $\bkb\rec \bka$.
Let $\bka_\dashv$ and
$\bka_i$, $1\le i\le m$  be defined as in (1.2) and (1.9), and  $k$ 
 the minimal index with $a_k=r$.
Then $k\ge2$ and $\bka_k= \bka$. By induction,
$v_{\bka_j}\bH$, $1\le j\le k-1$, are projective, and   
$v_{\boldkey a_\dashv}\bH$ is projective. Using the short exact sequence in 
(4.12), we have $v_{a_\dashv}\bH\cong\oplus_{j=1}^{k-1} v_{\bka_j}\bH\oplus V_k$
and hence, $V_{k}$ is projective. By (4.13), we have
$V_k\cong V_{k+1}\oplus v_{\bka_k}\bH$, and
consequently,
$v_{\bka}\bH=v_{\bka_k}\bH  $ is projective, proving (a).
Note from (4.7) and (4.12-3)  
that we have for any $r>0$ and $\bka\in\La[m,r]$
$$v_{\bka_\dashv}\bH\cong\oplus_{i=1}^m v_{\bka_i}\bH.\tag4.15$$

To see (b), we apply the tensor functor $-\otimes_{\bH_m^{r_1} } \bH$ 
to (4.15) and
 obtain, for $1\le r_1\le r$ and 
$\bka\in\La[m,r_1]$
$$v_{\bka_\dashv}\bH\cong\oplus_{j=1}^{m} v_{\bka_j}\bH,$$
from which (b) follows.

Putting $r_1=r$ in both (a) and (b), we see that
the $\bH$-module $\oplus_{\bka\in \La[m, r]}v_{\bka}\bH$ is a
projective generator. 
  By  standard results (see, e.g.,  \cite{CR, (3.54)}) 
the categories  of $\bH$-modules and $\varepsilon \bH \varepsilon$-modules
are Morita equivalent (cf. (3.11)). 
\qed\enddemo

\subhead 5. Some applications\endsubhead
We first recall the notion of multi-compositions and multi-partitions of $r$.
By definition,    $\bla=(\la^{(1)}, \cdots, \la^{(m)})$ is 
 called an $m$-composition (resp. $m$-partition) of $r$ if $\la^{(i)}$ is a 
composition (resp. partition)  for every $i$ with $1\le i\le m$ and 
$|\bla|=\sum_{i=1}^m |\la^{(i)}|=r$. 
Recall that  $\La(n,r)$ (resp. $\La(n, r)^+$) is  the set of all compositions
(resp. partitions) of $r$ with $n$ parts, and let,  for $m>0$ and $\mu
=(\mu_1,\cdots,\mu_m)\in\La(m,r)$,
$$\aligned
&\La(n,\mu)=\La (n_r, \mu_1)\times...\times\La (n_r, \mu_{m-1})\times
\La (n, \mu_m)\cr
&\La_m(n, r)=\bigcup_{\mu\in\La(m,r)}\La(n,\mu)\cr 
\endaligned 
$$ 
where $n_r$ is the maximum of $n$ and $r$. Define
$\La_m(n, r)^+$ similarly.
The set $\La_m(n,r)$ can be identified with $\La(N,r)$ by concatenation,
where $N=(m-1)n_r+n$. To distinguish them, for $\bla\in\La_m(n,r)$,
let $\bar\bla$ denote
the corresponding element in $\La(N,r)$.

For $\bla\in \La_m(n, \mu )$,  let $\boldkey a=[a_0, \cdots, a_m]$
be {\it the cumulative
norm sequence} of $\bla$, denoted cns$(\bla)$, where
$a_0=0$ and $a_i=\mu_1+\cdots+ \mu_i$ for $1\le i\le m$.
Let $e$ be the minimal integer $l$  such that 
$$
1+q+q^2 +\cdots +q^{l-1}=0.
$$
If such an integer $l$ does not exist, then set $e=\infty$.
A partition $\la$ is called $e$-regular if $\la$ has no 
non-zero part occurring $e$ or more times.
An $m$-partition $\bla=(\la^{(1)}, \cdots, \la^{(m)})$ is $e$-regular
if each $\la^{(i)}$, $1\le i\le m$,  is $e$-regular.

Let $\sH_F$ be the Hecke algebra of type $A_{r-1}$ over a field $F$, in which
$q$ is a primitive $e$-th root of 1.
Let $S^\la$ be the Specht module with respect to $\la$.
Then, by \cite{DJ1}, $S^\la$ has simple head 
 $D^\la$, if $\la$ is $e$-regular. (Note, for $e>r$, $D^\la=S^\la$.)
For a multi-partition $\bla=(\la^{(1)}, \cdots, \la^{(m)})$ of $r$
with $\cns(\bla)=\boldkey a$, 
let $S^{\la^{(1)}}\cdots S^{\la^{(m)}}$ ($\cong \otimes_{i=1}^m S^{\la^{(i)}}$)
be the corresponding Specht module for $\sH_F(\frak S_{\boldkey a})$,
and let $S^{\bla}$ be the right ideal
 $e_{\boldkey a} S^{\la^{(1)}}\cdots S^{\la^{(m)}} \bH_F$ of $\bH_F$.
By (4.14)(c), we know that  the $\bH$-module $S^{\bla}$ and 
$\varepsilon \bH \varepsilon$-module $S^{\bla}\varepsilon$ 
have isomorphic submodule
lattices. Since $S^{\bla}\varepsilon\cong S^{\la^{(1)}}\cdots S^{\la^{(m)}}$
and the latter has simple head if $\bla$ is $e$-regular,
so is $S^{\bla}$. Let $D^{\bla}$ denote the simple head of $S^{\bla}$.
The following is an immediate consequence of (4.14)(c)
(see \cite{DJ2, (5.3)} for the case $m=2$).

\proclaim{(5.1) Theorem} Let $F$ be a field and $f_{m, r}\neq 0$ in $F$.
Then the set 
$$
\{ D^{\bla}\mid \bla\in \La_m(r, r)^+ \   e\text{-regular} \}
$$
is a complete set of simple $\bH_F$-modules.
\endproclaim

The following result has been proved in  \cite{Ari}. Recall the Poincar\'e
polynomial $d_W$ introduced in (4.2)
 
 \proclaim{(5.2) Theorem} Let $\bH_F$ be the Ariki-Koike algebra 
 over a field $F$. Then $\bH_F$ is semisimple if and only if 
$d_W\neq 0$.
\endproclaim

\demo{Proof} Note that $d_W\neq 0$ is equivalent to $f_{m,r}\neq0$ and 
$e>r$.
By (4.3), we have $f_{m, r}\neq 0$ if and only if 
$f_{m, r, i}\neq 0$ for all $i$. 
If $f_{m, r}\neq 0$, then the categories of $\bH$-modules and
$\varepsilon\bH\varepsilon$-modules are Morita equivalent. By (3.10) 
and  \cite{DJ1, 4.3}, $\varepsilon\bH\varepsilon$ is semisimple if and only if
$e>r$, proving the ``if'' part.
Conversely, it suffices to look at the case $f_{m, r}=0$. Then
$f_{m, i, r}=0$ for some $i$, and
$z_{\boldkey r_i'}$ is not invertible. Therefore, 
$v_{\boldkey r_i'}\bH$ is not an idempotent ideal by (3.9).
Thus,
$\bH$ is not semisimple. \qed\enddemo

We finally look at a Morita theorem between  $q$-Schur$^m$
algebras and $q$-Schur algebras (i.e., the $q$-Schur$^1$ algebras). 
For any $\bla\in \La_m(n, \boldkey r)$  with $\cns(\bla) =\bka$, let 
$x_{\bla}=\pi_\bka x_{\bar\bla}$.

\proclaim{(5.3) Lemma}
Let $R$ be an integral domain and assume that $f_{m, r}\in R$ is a unit.
For any $\bla\in \La_m(n, r)$ with $\cns(\bla)=\bka$, we have
$e_{\bka} x_\bla\bH=e_{\bka }x_{\bar\bla}\bH$. \endproclaim 

\demo{Proof} Since $\pi_\bka\bH  \subseteq \bH$, we have $e_{\bka}x_\bla\bH\subseteq e_{\bka}x_{\bar\bla}\bH$.
By (1.8), (3.1)(a) and (3.6)-(3.7), we have 
$e_{\boldkey a} x_{\bar\bla}  =v_{\boldkey a} 
h_{\boldkey a'} z_{\boldkey a}^{-1} x_{\bar\bla}
=v_{\boldkey a} h_{\boldkey a'}x_{\bar \bla} z_{\boldkey a}^{-1} 
=v_{\boldkey a} h h_{\boldkey a'}  z_{\boldkey a}^{-1}
=
x_{\bar \bla} v_{\boldkey a}  h_{\boldkey a'} z_{\boldkey a}^{-1} 
=x_{\bar \bla}   e_{\boldkey a}
$,
where $h\in \sH(\frak S_{\bka'})$ given by
 $ h h_{\bka'}=h_{\bka'} x_{\bar\bla}$.
Thus, $e_\bka x_{\bla}\bH=x_{\bar\bla} e_{\bka} \pi_{\bka}\bH=x_{\bar\bla} v_{\bka} h_{\bold a'} \pi_{\bka}\bH  
\supseteq x_{\bar\bla} v_{\bka} h_{\bold a'} 
v_\bka \bH =x_{\bar\bla} v_{\bka} z_{\bold a'}\bH=x_{\bar\bla} e_{\bka}\bH$, 
proving (5.3).\qed\enddemo

 Some special case for $m=2$ of the second part of the following result
has been 
discussed by Gruber and Hiss \cite{GH} in the context of representations
finite groups of Lie type.

\proclaim{(5.4) Corollary} Let $\bH$ be the Ariki-Koike algebra over an 
integral domain  $R$, in which $f_{m, r}$ is a unit. 
Let $\bka=\cns(\bla)$ for $\bla\in \La_m(n, r)$. Then 
$$\aligned
\End_{\bH}\left(\oplus_{\bla\in\La_m(n,r)} e_{\bka}
x_{\bla}  \bH\right)
&=\End_{\bH}\left(\oplus_{\bla\in\La_m(n,r)} e_{\bka}
x_{\bar\bla}\bH\right)\cr
&\cong
\oplus_{\mu\in\La(m,r)}\End_{\sH(\fS_\mu)}\left(\oplus_{\bla\in\La(n,\mu)}
x_{\bar{\bla} }  \sH(\fS_\mu)\right).\cr
\endaligned$$ 
\endproclaim

\demo{Proof} The first equality follows from (5.3).
Using standard results
(see \cite{AF, (21.2)} or \cite{DPS, (0.1)}), 
we know that 
$$\End_{\bH}\left(\oplus_{\bla\in\La_m(n,r)} e_{\bka}
x_{\bar\bla}\bH\right)\cong\End_{\varepsilon \bH\varepsilon}
\left(\oplus_{\bla\in\La_m(n,r)} e_{\bka}
x_{\bar\bla}\bH\varepsilon\right).$$
For $\bla=(\la^{(1)}, \cdots, \la^{(m)})\in\La_m(n,r)$,
write $\cns(\bla)=\boldkey a=[a_i]$. Then we have 
$x_{\bar{\bla}}\sH(\frak S_{\boldkey a})\cong
x_{\la^{(1)}} \sH(\frak S_{a_1})\otimes_R \cdots\otimes_R
x_{\la^{(m)}} \sH(\frak S_{a_m-a_{m-1}})$,
and by the proof of (5.3), 
$e_\bka x_{\bar\bla}=x_{\bar\bla} e_{\bka}$. Therefore,
$e_{\bka} x_{\bar{\bla} }  \bH \varepsilon
=x_{\bar \bla}e_{\bka}\sH( \frak S_{\boldkey a})\cong x_{\bar\bla}\sH(\fS_{\bka})$.
By (3.10), we have
$$\aligned 
\End_{\varepsilon \bH\varepsilon  }\left(\oplus_{\bla\in\La_m(n,r)}e_{\bka}
x_{\bar{\bla} }  \bH \varepsilon\right)
&\cong\oplus_{\bka\in\La[m,r]}\End_{e_\bka\bH e_\bka}
\left(\oplus_{\bla\in\La(n,\Theta(\bka))}
x_{\bar{\bla} } 
e_{\bka} \bH e_{\bka}\right)\cr
&\cong\oplus_{\mu\in\La(m,r)}\End_{\sH(\fS_\mu)}\left(\oplus_{\bla\in\La(n,\mu)}
x_{\bar{\bla} }  \sH(\fS_\mu)\right).\cr
\endaligned
$$
Here $\Theta$ is defined in (1.1).
\qed\enddemo

The endomorphism algebra 
$\bS_R^m(n,r)=\End_{\bH}\left(\oplus_{\bla\in\La_m(n,r)}
x_{\bla}  \bH\right)$ is called a $q$-Schur$^m$ algebra in \cite{DR2},
which is a 
cyclotomic $q$-Schur algebra in the sense of \cite{DJM}.
It is proved in \cite{DR2} that a Borel type subalgebra of 
a $q$-Schur$^m$ algebra
is isomorphic to a Borel subalgebra of a $q$-Schur algebra.
Now, with the invertibility of the polynomial $f_{m,r}$, we establish
below a Morita equivalence  between the module categories
of a $q$-Schur$^m$ algebra and a direct sum
of tensor products of certain $q$-Schur algebras. 

\proclaim{(5.5) Theorem} 
Let $\bH=\bH_\sO$ be the Ariki-Koike algebra over
a discrete valuation ring $\sO$. Then the $q$-Schur$^m$ algebra $\bS_\sO^m(n,r)$
is Morita equivalent to the algebra
$$
\bigoplus_{\mu\in\La(m,r)}
\bS_\sO^1(n_r,\mu_1)\otimes...\otimes\bS_\sO^1(n_r,\mu_{m-1})\otimes
\bS_\sO^1(n,\mu_m).
$$
\endproclaim

\demo{Proof} 
We first note the isomorphism
$$\bigoplus_{\mu}\End_{\sH_\sO
(\fS_\mu)}\left(\oplus_{\bla\in\La(n,\mu)}
x_{\bar{\bla} }  \sH_\sO(\fS_\mu)\right)\cong 
\bigoplus_{\mu} \otimes_{i=1}^{m-1} 
\bS_\sO^1(n_r,\mu_i) \otimes \bS_\sO^1(n, \mu_m),\tag5.6
$$
where $\mu\in\La(m,r)$.
The $q$-Schur algebra of 
bidegree $(n, r)$ is a
(integral) quasi-hereditary algebra, 
whose  simple modules are parametrized by $\La(n, r)^+$
(see \cite{DS, \S2}). Therefore, by \cite{Wi, (1.3)} (or
a cellular basis argument \cite{GL},\cite{DR1}), the algebra  
$$
\bigoplus_{\mu\in\La(m,r)} \otimes_{i=1}^{m-1} 
\bS_\sO^1(n_r,\mu_i) \otimes \bS_\sO^1(n, \mu_m)
$$
 is a quasi-hereditary algebra,
 whose simple modules are indexed by $\La_m(n, r)^+$.
Thus, PIMs are indexed by  $\La_m(n, r)^+$, too. 
Using Fitting's Lemma,  
the non-isomorphic indecomposable direct summands (The existence
of these modules follows from Heller's result, \cite{CR, (30.18iii)}.)
of 
$\oplus_{\bla\in\La_m(n,r)}  x_{\bar \bla } e_{\bka} \bH_\sO e_{\bka}$ are 
indexed by 
$\La_m(n, r)^+$. 
Therefore, the non-isomorphic indecomposable direct summands of 
$\oplus_{\bla\in\La_m(n,r)}  e_{\bka}x_{\bla }\bH_\sO$ are 
indexed by 
$\La_m(n, r)^+$ by (5.4).
Since 
$e_\bka x_{\bar\bla}=x_{\bar\bla} e_{\bka}$ and $e_\bka^2=e_\bka$, we have 
$e_\bka x_{\bla} \bH_\sO\oplus (1-e_\bka) x_{\bla} \bH_\sO=x_\bla\bH_\sO$. 
It follows that every direct summand of $e_\bka x_{\bla} \bH_\sO$
is a direct summand of $x_\bla\bH_\sO$.
Now, by the quasi-heredity of  the $q$-Schur$^m$ algebra
(see \cite{DR2, (5.10)}),
 the non-isomorphic indecomposable direct summands of 
$\oplus_{\bla\in\La_m(n,r)}  x_{ \bla } \bH_\sO$ are 
indexed by 
$\La_m(n, r)^+$.
Therefore, both $\oplus_{\bla\in\La_m(n,r)}  x_{\bla } \bH_\sO$
and $\oplus_{\bla\in\La_m(n,r)}  e_\bka x_{\bla } \bH_\sO$
have the same non-isomorphic indecomposable direct summands.
Consequently, the $q$-Schur$^m$ algebra $\bS_\sO^m(n,r)$ 
is Morita equivalent to
$\End_{\bH_\sO}\left(\oplus_{\bla\in\La_m(n,r)} e_{\bka} x_{\bla} \bH_\sO\right)$.
Now, the required Morita equivalence follows from
(5.4) and (5.6).
\qed\enddemo

\end